\documentclass[11pt, leqno]{article}
\usepackage{graphicx}
\usepackage{textcomp}
\usepackage{amsmath,amsfonts,amscd,amsthm,amssymb}
\usepackage{epsfig}
\usepackage{makeidx}
\usepackage{color} 
\usepackage{subfig}

\makeindex
\newcommand{\ncom}{\newcommand}
\newcommand{\be}{\begin{eqnarray*}}
\newcommand{\ee}{\end{eqnarray*}}

\newcommand{\ben}{\begin{eqnarray}}
\newcommand{\een}{\end{eqnarray}}

\newcommand{\tl}{\tilde}

\newcommand{\pr}{\partial}
\newcommand{\nab}{\nabla}

\newcommand{\bu}{{\bf u}}
\newcommand{\bpsi}{{\bf \psi}}
\newcommand{\bU}{{\bf U}}

\newcommand{\bw}{{\bf w}}
\newcommand{\e}{{\bf e}}
\newcommand{\bs}{{\bf \sigma}}

\newcommand{\bJ}{{\bf J}}

\newcommand{\bv}{{\bf v}}

\newcommand{\bta}{\mbox{\boldmath $\eta$}}
\newcommand{\bzeta}{\mbox{\boldmath $\zeta$}}

\newcommand{\bphi}{\mbox{\boldmath $\phi$}}
\newcommand{\brho}{\mbox{\boldmath $\rho$}}

\newcommand{\bxi}{\mbox{\boldmath $\xi$}}
\newcommand{\bH}{{\bf {H}}}

\newcommand{\bL}{\bf {L}}

\newcommand{\f}{{\bf {f}}}
\newcommand{\bvh}{{\bf v_h}}

\newtheorem{tdf}{Theorem}[section]

\newtheorem{ldf}{Lemma}[section]
 
\ncom{\ul}{\underline}
\ncom{\beq}{\begin{equation}}
\ncom{\eeq}{\end{equation}}
\ncom{\bea}{\begin{eqnarray*}}
\ncom{\eea}{\end{eqnarray*}}
\ncom{\beqa}{\begin{eqnarray}}
\ncom{\eeqa}{\end{eqnarray}}
\ncom{\nno}{\nonumber}
\ncom{\non}{\nonumber}
\ncom{\ds}{\displaystyle}
\ncom{\half}{\frac{1}{2}}
\ncom{\mbx}{\makebox{.25cm}}
\ncom{\hs}{\mbox{\hspace{.25cm}}}
\ncom{\rar}{\rightarrow}
\ncom{\Rar}{\Rightarrow}
\ncom{\noin}{\noindent}
\ncom{\sz}{\scriptsize}
\ncom{\Sgm}{\Sigma}
\ncom{\psgm}{\sigma^{\prime}}
\ncom{\dt}{\delta}
\ncom{\Dt}{\Delta}
\ncom{\lmd}{\lambda}
\ncom{\Lmd}{\Lambda}
\ncom{\Th}{\Theta}
\ncom{\eps}{\epsilon}
\ncom{\pcc}{\stackrel{P}{>}}
\ncom{\lp}{\stackrel{L_{p}}{>}}
\ncom{\sspan}{{\rm\,span}}
\ncom{\re}{{\rm Re\,}}
\ncom{\im}{{\rm Im\,}}
\ncom{\sgn}{{\rm sgn\,}}
\ncom{\ba}{\begin{array}}
\ncom{\ea}{\end{array}}
\ncom{\integ}[4]{\int_{#1}^{#2}\,{#3}\,d{#4}}
\ncom{\vspan}[1]{{{\rm\,span}\{ #1 \}}}
\ncom{\dm}[1]{ {\displaystyle{#1} } }
\ncom{\ri}[1]{{#1} \index{#1}}

\newtheorem{remark}{\bf Remark}[section]

\newtheorem{example}{Example}[section]

\newtheoremstyle
    {remarkstyle}
    {}
    {11pt}
    {}
    {}
    {\bfseries}
    {:}
    {     }
    {\thmname{#1} \thmnumber{#2} }

\theoremstyle{remarkstyle}

\def \R{{{\rm I{\!}\rm R}}}

\begin{document}
\title
 {A priori error estimates of fully discrete finite element Galerkin method for 
 Kelvin-Voigt viscoelastic fluid flow model}
 \author
{Saumya Bajpai \\
 School of Mathematics and Computer Science, \\
Indian Institute of Technology Goa, Ponda-403401, India\\
and\\
Ambit K. Pany\\
Department of Mathematics, Gandhi Institute for Technological Advancement,\\
Bhubaneswar-752054, India
}
\maketitle 
\begin{abstract}
In this article, a finite element Galerkin method is applied to the Kelvin-Voigt viscoelastic 
fluid model, when its forcing function is in $L^{\infty}(\bL^2)$. 
Some new {\it a priori} bounds for the velocity as well as for the pressure are derived which 
are independent of inverse powers of the retardation time $\kappa$. Optimal error estimates for the velocity 
in $L^{\infty} (\bL^2)$ 
as well as in $L^{\infty}(\bH^1_0)$-norms and for the pressure in $L^{\infty}(L^2)$-norm of the semidiscrete method 
are discussed 
which hold uniformly with respect to $\kappa$ as $\kappa\rightarrow 0$ with the initial condition only in 
$\bH^2\cap\bH_0^1$. 
Further, under uniqueness 
condition, 
these estimates are shown to be uniformly in time as $t \mapsto \infty$. For the complete 
discretization of the semidiscrete 
system, a first-order 
accurate backward Euler method is applied and 
fully discrete optimal error estimates are established. Finally, numerical experiments are conducted to 
verify the theoretical results. The results derived in this article are sharper 
than those derived earlier for finite element analysis of the Kelvin-Voigt fluid model in the sense that the error estimates 
in this article hold true uniformly even as $\kappa\rightarrow 0$. 
\end{abstract}

\section{\normalsize\bf Introduction}
 
The equations of motion arising from the Kelvin-Voigt model give rise to the following system 
of partial differential equations :
\begin{eqnarray}
\label{1.1}
\frac {\partial \bu}{\partial t}+ \bu\cdot \nabla \bu - \kappa \Delta \bu_t 
-\nu \Delta \bu  + \nabla p
=\f(x,t),\,\,\, x\in \Omega ,\,\,t>0,
\end{eqnarray}
and incompressibility condition
\begin{eqnarray}
\label{1.2}
\nabla \cdot \bu=0,\,\,\,x\in \Omega,\,t>0,
\end{eqnarray}
with initial and boundary conditions
\begin{eqnarray}
\label{1.3}
\bu(x,0)= \bu_0 \;\;\;\mbox {in}\;\Omega,\;\;\;\;\; \bu=0,\;\; \;
\mbox {on}\; \partial \Omega,\; t\ge 0.
\end{eqnarray}
Here, $\Omega$ is a bounded convex polygonal or polyhedral domain in $\R^d, d=2,3$  with 
boundary $\partial \Omega$, $\bu=(u_1,u_2)$ 
$(\mbox{or}~ \bu=(u_1,u_2,u_3))$ represents the velocity vector, $p$ is 
the pressure of the fluid, $\f$ is the external force, $\nu>0$ denotes the kinematic coefficient of viscosity and 
$\kappa$ is the retardation time. For a more physical description and applications 
of the model, one may refer \cite{bs05}-\cite{css02}, \cite{AP1} and literature therein. Based on the the proof 
techniques of Ladyzenskaya \cite{L69}, Oskolkov and his collaborators \cite{AP1}, \cite{AP2}, \cite{OR1}, \cite{OR2} have 
discussed the existence of a unique global ``almost `` classical solution for the initial and boundary value problem 
(\ref{1.1})-(\ref{1.3}) for various 
assumptions on the right-hand side function $\f$ and for all time $t>0$.\\  
 \noindent
There is a considerable amount of literature devoted to the numerical approximations of Kelvin-Voigt fluid flow model, 
see \cite{bnpdy}-\cite{BN},  \cite{KBP}, \cite{AP3}, \cite{PDY}-\cite{AKP}. 
In \cite{AP3}, Oskolkov  has applied the spectral 
Galerkin approximation to the problem (\ref{1.1})-(\ref{1.3}) and has proved the convergence for $t\geq 0$ 
with the assumption that the solution is asymptotically stable as $t\rightarrow\infty$. Further, the author established 
optimal error estimates in $L^{\infty}({\bH_0^1})$-norm, which are local in time, since the 
constants appearing in error bounds involve exponential in time terms. Later on, as an improvement to the Oskolkov work, 
Pani et al. \cite{PDY} have 
established $L^{\infty}(\bL^2)$ and $L^{\infty}(\bH_0^1)$-norms optimal error estimates for the spectral 
Galerkin method applied to 
(\ref{1.1})-(\ref{1.3}), which are valid uniformly in time under uniqueness assumption. They further applied modified 
nonlinear Galerkin method to (\ref{1.1})-(\ref{1.3}), and have established optimal uniform in time {\it a priori} 
error estimates with the assumption of uniqueness condition. They have also observed 
the superconvergence phenomenon in $L^{\infty}(\bH_0^1)$-norm for both spectral Galerkin method and 
modified nonlinear spectral 
Galerkin method. Note that, the constants appearing in error estimates derived in
\cite{bnpdy}-\cite{BN},  \cite{KBP}, \cite{AP3}, \cite{PDY} depend on $\kappa^{-r}$, 
for $r\geq 2$ which may blow as $\kappa\rightarrow 0$. \\
\noindent 
 In \cite{ANS6}, the authors have applied semidiscrete finite element Galerkin method to the problem 
(\ref{1.1})-(\ref{1.3}) and have established some new uniform in time a priori bounds for 
the weak solution. It can be observed that the constants appearing in {\it a priori} bounds for the weak solution 
are independent of inverse powers of $\kappa$ 
which is an improvement over the results 
derived in earlier articles 
related to the regularity estimates for the weak solution of this model. 
Further, using these {\it a priori} estimates, they have established optimal error estimates for 
the velocity in 
$L^{\infty}(\bL^2)$ as well as in $L^{\infty}(\bH_0^1)$- norms and 
for the pressure in $L^{\infty}(L^2)$-norm, when the forcing function $\f\in L^{\infty}(\bL^2)$. Here, it can be noted that 
they have achieved an improvement in the error estimates in powers of $\kappa$ as the constants in error bounds depend only on $\kappa^{-1/2}$.\\
\noindent
 As an extension to the work in \cite{ANS6}, 
Pany et al. \cite{P4}, \cite{AKP} have employed a linearized first order backward Euler method 
and a second order backward difference scheme 
for the time discretization of the problem (\ref{1.1})-(\ref{1.3}) with $\f\in L^{\infty}(\bL^2)$ and 
have derived {\it a priori} bounds for the discrete solution in the Dirichlet norm 
using a combination of discrete Gronwall's 
lemma and Stolz-Cesaro's classical result for sequences. Then, making use of these {\it a priori} estimates for the solution, 
they have established fully discrete optimal error estimates for the 
velocity and pressure, which hold true uniformly in time under uniqueness assumption. In 
\cite{AKP}, the author has also mentioned that assuming the solution is smooth enough, that is, $\bu_0\in \bH^3\cap {\bH_0^1}$ 
with $\Delta \bu_0=0$ on 
$\partial\Omega$, the optimal error estimates 
independent of $\kappa$ can be achieved following the similar analysis as in \cite{ANS6}-\cite{AKP}. 
For the articles related to the finite element analysis 
of the problem (\ref{1.1})-(\ref{1.3}) 
with the right-hand side forcing function $\f=0$, one may refer to \cite{bnpdy}-\cite{BN}. For the papers 
containing the similar results for the Navier-Stokes and Oldroyd models, see 
 \cite{AO}, \cite{GR}-\cite{YH}, \cite{PY}, \cite{PYD}, \cite{WLH1}, \cite{WLH2} and literature, referred therein.\\
\noindent
Since the Kelvin-Voigt fluid is characterized by the fact that after instantaneous removal of the stresses, 
the velocity of the fluid does not vanish instantaneously but dies out like $\exp(-\kappa^{-1}t)$ \cite{AP2}, it is worthwhile to discuss the behavior of the solution 
as $\kappa\rightarrow 0$ and as $t\rightarrow \infty$. Moreover, this model can 
be thought of as a $\kappa$ regularization of the Navier-Stokes model (\cite{KLR}, \cite{L69}). 
Based on these observations, in this article, we mainly aim 
at recovering optimal error estimates which are valid uniformly 
in time as well as in 
retardation time $\kappa$ 
under realistically assumed minimum regularity assumption on the exact  solution with $\bu_0\in\bH^2\cap \bH_0^1$ and $\f$, $\f_t \in L^{\infty}(\bL^2)$. 

\noindent
The main contributions of the present article are as follows:\\
\noindent
(i) Some new regularity results for the higher order time derivatives of the weak solution are derived which are valid 
uniformly in time. Further, these estimates are shown to be uniformly in $\kappa$ as 
$\kappa \rightarrow 0$ under minimum regularity assumptions $\bu_0\in\bH^2\cap \bH_0^1$ and $\f$, 
$\f_t \in L^{\infty}(\bL^2)$. Here, it can be noted that the introduction of weight 
function $\bs(t)=\min\{1,\,t\}e^{2\alpha t}$ plays a key role in handling 
the regularity issues at $t=0$.  \\
(ii) Using the Sobolev-Stokes projection defined earlier in \cite{bnpdy}, 
fully discrete optimal error estimates in $L^{\infty}(\bL^2 )$ and $L^{\infty}(\bH_0^1 )$-norms for 
the finite element velocity approximation and in $L^{\infty}(L^2 )$-norm for the 
finite element pressure approximation are established.  It is further proved that these error estimates 
hold uniformly as $\kappa\rightarrow 0$. Here, we would like to highlight an important point that we have 
resorted to a simple observation (mentioned in Remarks \ref{lr1}, \ref{lr2}) in order to derive the estimates 
involving weight function $\bs(t)$ which 
plays an important role in achieving uniform estimates in terms of $\kappa$.
\\
(iii) Since the error bounds derived in (ii) involve exponential in time terms, it is 
further established that under the assumption of uniqueness condition, the error estimates are uniformly in time.\\
(iv) Numerical results are presented to validate our theoretical findings. 
Moreover, it is depicted that the order of convergence does not degenerate as $\kappa\rightarrow 0$ confirming the 
results in $(ii).$\\

\noindent
Note that, the results in this article are substantial improvements over the results available in literature 
related to the finite element error analysis of the Kelvin-Voigt model in the sense that we are able to establish 
error bounds which do not involve inverse powers of $\kappa$. As a consequence, the error 
estimates do not blow up as $\kappa\rightarrow 0$. 
The main difficulty in making error estimates independent 
of $\kappa$ arises due to the lack of regularity of solution at $t=0$. In order to overcome this difficulty, 
we introduce various powers of weight function $\bs(t)$ which takes care of  regularity issues of the
solution at $t=0$.\\
\noindent
The remaining part of the article consists of the following sections. In Section {\bf 2}, some preliminaries to be 
used in the subsequent sections are introduced and some new regularity results for the weak solution are derived. 
In Section {\bf 3}, assumptions on finite element spaces 
to determine the discrete solution are presented and semidiscrete finite element approximations are defined.  The main results of the article are also stated.
Section {\bf 4} deals with the optimal error estimates for velocity and pressure. In 
Section {\bf 5}, full discretization is achieved by using the backward Euler method. 
Section {\bf 6} presents some numerical results which confirm our theoretical findings. 
Finally, Section {\bf 7} concludes the article by briefly summarizing the results.



\section{\normalsize \bf Preliminaries and Weak formulation}
\setcounter{tdf}{0}
\setcounter{ldf}{0}
\setcounter{cdf}{0}
\setcounter{equation}{0}
We denote $\mathbb {R}^d$ $(d=2, 3)$-valued function spaces using bold 
face letters, that is, ${\bf H}_0^1 = (H_0^1(\Omega))^d$, ${\bf L}^2 = (L^2(\Omega))^d$ and ${\bf H}^m=(H^m(\Omega))^d,$ 
where ${\bL}^2(\Omega)$ is the space of square integrable functions defined in $\Omega$ with inner product $(\bphi,\bpsi) = \displaystyle{\int_0^t}\bphi({\bf x})\bpsi({\bf x})dx$ 
and norm $\|\bphi\| =\left(\displaystyle{\int_0^t} |\bphi({\bf x})|^2 dx\right)^{1/2}$. 
Further, $\bH^m(\Omega)$ denotes the standard Hilbert Sobolev space of order $m \in N^{+}$ 
with norm $\|\bphi\|_m = \displaystyle{\sum_{|\alpha|\leq m}}\left (\displaystyle{\int_0^t}|D^{\alpha} \bphi|^2 dx\right)^{1/2}$. The space $\bH_0^1$ is equipped with a norm $\|\nabla\bv\|= \left({\sum_{i,j=1}^{d}}
(\partial_j v_i, \partial_j v_i)\right)^{1/2}=\left({\sum_{i=1}^{d}}
(\nabla v_i, \nabla v_i)\right)^{1/2}$.
Given a Banach space $X$ endowed with norm ${\parallel \cdot \parallel}_X$, 
let $L^p(0,T;X)$ be the space of all strongly measurable functions $\phi:(0,T) \rightarrow X$ 
satisfying $\displaystyle{\int_0^{T}} {\parallel \phi(s) \parallel}^p_X\,ds< \infty$ if 
$ 1\leq p<\infty$ and 
for $p=\infty$, $\rm{\displaystyle{\mathop{ess \sup}_{t \in(0,T)}} {\parallel \bphi(t) \parallel}_X< \infty}$. 
Also, we define the divergence free spaces
 \begin{eqnarray}
 {\bf J}&=& \{\bphi \in {\bf {L}}^2 :\nabla \cdot \bphi = 0\;\;
{\mbox {\rm in}}
 \;\; \Omega, \;\;\bphi \cdot {\bf {n}} |_{\partial \Omega} = 0\;\;
 {\mbox {\rm holds} \;  {\rm weakly}} \},\nonumber\\
{\bf J}_1 &=& \{{\bphi} \in {\bf {H}}_0^1 : \nabla \cdot \bphi = 0\},\nonumber
 \end{eqnarray}
where ${\bf {n}}$ is the unit outward normal to the boundary
$\pr \Omega$ and $\bphi \cdot {\bf {n}} |_{\pr \Omega} = 0$ should
be understood in the sense of trace in $\bH^{-1/2}(\partial \Omega)$,
see \cite{temam}. Let $H^m/{\rm I\!R} $ be the quotient space with norm $\| \phi\|_{H^m /{\rm I\!R}}
 = \inf_{c\in{\rm I\!R} }\| \phi+c\|_m$. For $m=0$, it is denoted by $L^2 /{\rm I\!R}$. Now, define 
 $P : \bL^2 \rightarrow \bJ$ as the $\bL^2$-orthogonal projection.\\
\noindent
Throughout this article, we make the following assumptions: \\
\noindent
({\bf A1}). 
Setting
$-\tl{\Delta} = -P\Delta:{\bf J}_1 \cap {\bf H}^2
\subset {\bf J} \rightarrow  {\bf J}
$
as the Stokes operator, assume that the following regularity result holds:
\ben \label{2.1a} 
\| \bv\|_2 \le C \| \tl{\Delta} \bv \| \;\;\;\;\forall \bv\in {\bf J}_1
\cap {\bf H}^2.
\een
The above assumption is valid as the domain $\Omega$ is a convex polygon or convex polyhedron.
It can be noted that the following Poincar\'e inequality \cite{HR82} holds true:
\ben \label{2.1*} 
\|\bv\|^2 \le \lambda_1^{-1} \| \nabla \bv \|^2
\;\; \forall \bv \in {\bf H}_0^1(\Omega),
\een
where $ \lambda_1^{-1}$, is the best possible positive constant depending on the domain $\Omega.$
Further, observe (see, \cite{HR82}) that 
\ben \label{2.1**} 
\| \nabla \bv\|^2 \le \lambda_1^{-1} \| \tl{\Delta} \bv\|^2\;\;
\forall \bv  \in {\bf J}_1 \cap {\bf H}^2.
\een
\noindent
({\bf A2}). There exists a positive constant $ M $ such that the initial velocity $\bu_0$ and the external 
force $\f$ satisfy for 
$t \in (0,\,T]$ with $0<T<\infty$
\be
\bu_0 \in {\bf H}^2\cap {\bf  J}_1,\, \f,\;\f_t \in L^{\infty}(0,\,T;\,\bL^2)\,\, \rm{with}
\,\, \|\bu_0\|_2 \le M,\,\, 
\rm{\displaystyle{\mathop{ess\,\sup}_{0<t\leq T}}}\,\{\|\f(\cdot,t)\|,\|\f_t(\cdot,t)\|\}\leq M.\nonumber
\ee
\noindent
The weak formulation of (\ref{1.1})-(\ref{1.3}) is to find $(
\bu(t), p(t))\in \bH_0^1\times L^2 /{\rm I\!R}$, such that $\bu(0)= \bu_0$ and for $t>0$ 
\begin{eqnarray} \label{2.2atw}
 \left.
 \begin{array}{rcl}
&&(\bu_t, \bphi)+\kappa\, (\nabla \bu_t, \nabla\bphi) +\nu\, (\nabla \bu, \nabla\bphi)+( \bu\cdot\nabla\bu, \bphi)-( p,\nabla \cdot \bphi) 
=({\bf f},\bphi)
\;\;\;\;\; \forall \bphi \in {\bf H}_0^1,\\
&&(\nabla \cdot \bu,\chi)=0\;\;\;\forall \chi\in L^2.
 \end{array}
 \right\}
 \end{eqnarray}
\noindent
Equivalently, find  $\bu(t) \in {\bf J}_1 $ such that for $\bu(0)= \bu_0$, $t>0$,
\begin{align}
\label{2.3atw}
(\bu_t, \bphi) +\kappa\, (\nabla \bu_t, \nabla\bphi)+\nu\,(\nabla \bu, \nabla\bphi)+( \bu\cdot\nabla\bu, \bphi)= ({\bf f},\bphi)\;\;\;\forall \bphi \in {\bf J}_1.
 \end{align}
  For
$\bv, \bw, \bphi \in \bH_0^1$, define the bilinear form $a(\cdot,\cdot)$ as
$$
a(\bv, \bphi) = (\nabla \bv, \nabla \bphi)
$$
and the trilinear form $b(\cdot,\cdot,\cdot)$ as 
$$
b(\bv, \bw,\bphi)= \frac{1}{2} (\bv \cdot \nabla \bw , \bphi)
- \frac{1}{2} (\bv \cdot \nabla \bphi, \bw).
$$ 
Note that,  for $\bv \in \bJ_1$, $\bw$, $\bphi \in \bH_0^1$, $b(\bv, \bw, \bphi)=( \bv \cdot \nabla \bw, \bphi)$. 
Because of antisymmetric property of the trilinear form, it is easy
to verify that
$$b (\bv, \bw, \bw )= 0\,\,\, \forall \bv, \bw \in \bJ_1.$$
We present below in Lemma \ref{l1}, some {\it a priori} bounds for the weak solution pair $(\bu,\,p)$ which will be used in our subsequent error analysis. Since the estimates in (\ref{1r1}) and (\ref{1r2}) are already derived in \cite{ANS6}, we only provide proof of (\ref{1r3}).
\begin{ldf}\label{l1}[\cite{ANS6}, pp 241, 244]
Let the assumptions {\rm{({\bf {A1}})}}-{\rm{({\bf {A2}})}} hold. Then, there exists a 
positive constant $C =C(\nu,\alpha, \lambda_1, M)$ such that 
for $\displaystyle{0\leq\alpha< \frac{\nu \lambda_1}{4\big(1+\lambda_1\kappa\big)}}$ the following 
estimates hold true:
\begin{align}
& \displaystyle{\sup_{0<t< \infty}} \{\|\bu(t)\|_2+\|\bu_t(t)\|+\kappa \,\|\tilde\Delta\bu_t(t)\|
+\|p(t)\|_{\bH^1/{\mathbb R}}\}\leq C,\label{1r1}\\
&  e^{-2\alpha t}\int_0^t e^{2\alpha s}\|\bu_s(s)\|_1^2\;ds\leq C,\label{1r2}\\
&\bs^{-1}(t)\displaystyle{\int_0^t}e^{2\alpha s}
(\| \bu(s)\|_2^2+\|p(s)\|_{\bH^1/{\mathbb R}}^2+\kappa\,\|\tilde\Delta\bu_s(s)\|^2)\;ds\leq C,\label{1r3}
\end{align}
where $\tau(t):=\min\{t,\,1\}$ and $\bs(t):=\tau(t)e^{2\alpha t}$. 
\end{ldf}
\noindent{\it Proof.} We know that 
$\bs(t)=\tau(t)\,e^{2\alpha t}$, where $\tau(t)=\min\{t,1\}= t\,\text{or}\,1$. Hence, $\bs(t)=e^{2\alpha t}$ or 
$\bs(t)=t\,e^{2\alpha t}$. \\
Now, consider the following two cases:  \\
\noindent
{\bf Case 1:} $\bs(t)=e^{2\alpha t}$.
Then, 
 \begin{align}\label{n1}
 \frac{1}{\bs}\displaystyle{\int_0^t}e^{2\alpha s} \|\tilde\Delta\bu(s)\|^2ds&\leq   \frac{1}{e^{2\alpha t}} 
 \displaystyle{\sup_{0<t< \infty}}\|\tilde\Delta\bu(t)\|^2\displaystyle{\int_0^t}e^{2\alpha s} ds.
 \end{align}
 A use of (\ref{1r1}) in (\ref{n1}) leads to
 \begin{align}
 \frac{1}{\bs}\displaystyle{\int_0^t}e^{2\alpha s} \|\tilde\Delta\bu(s)\|^2ds & \leq C e^{-2\alpha t}
 \bigg( \frac{e^{2\alpha t}-1}{2\alpha}\bigg)\leq C(\alpha)(1-e^{-2\alpha t})\leq C. 
 \end{align}
{\bf Case 2:} $\bs(t)=t\,e^{2\alpha t}$. Again use (\ref{1r1}) and well known facts of series to obtain
\begin{align}
\frac{1}{\bs}\displaystyle{\int_0^t}e^{2\alpha s} \|\tilde\Delta\bu(s)\|^2ds&\leq   \frac{1}{t e^{2\alpha t}} 
 \displaystyle{\sup_{0<t< \infty}}\|\tilde\Delta\bu(t)\|^2\displaystyle{\int_0^t}e^{2\alpha s} ds\nonumber\\
 & \leq \frac{C}{t e^{2\alpha t}}\bigg( \frac{e^{2\alpha t}-1}{2\alpha}\bigg)\nonumber\\
 &=\frac{C}{2\alpha t e^{2\alpha t}}\bigg(\frac{2\alpha t}{1!}+\frac{(2\alpha t)^2}{2!}
 +\frac{(2\alpha t)^3}{3!}+\cdots\bigg)\nonumber\\
&\leq \frac{C}{ e^{2\alpha t}}\,
 \bigg(1+\frac{2\alpha t}{1!}+\frac{(2\alpha t)^2}{2!}+\cdots\bigg)
 \leq \frac{C}{ e^{2\alpha t}}\,e^{2\alpha t}\leq C.
 \end{align}
Therefore, considering the above two cases, we arrive at 
\begin{align}
 \bs^{-1}(t)\displaystyle{\int_0^t}e^{2\alpha s}\| \bu(s)\|_2^2 ds\leq C.\nonumber
\end{align}
\noindent
Following the similar sets of arguments as above, we obtain 
\begin{align}
 \bs^{-1}(t)\displaystyle{\int_0^t}e^{2\alpha s}(\|p(s)\|_{\bH^1/{\mathbb R}}^2+\kappa\,\|\tilde\Delta\bu_s(s)\|^2) 
 ds\leq C\nonumber
\end{align}
and this completes the remaining part of the proof.\hfill{$\Box$}\\
\noindent
In the next lemma, we derive {\it a priori} bounds for the highest order time derivatives of weak solution for the problem 
(\ref{2.2atw}).
\begin{ldf}\label{l2}
  Let the assumptions {\rm{({\bf {A1}})}}-{\rm{({\bf {A2}})}}  hold. Then, there exists a 
positive constant $C =C(\kappa,\nu,\alpha, \lambda_1, M)$ such that 
for $\displaystyle{0\leq\alpha< \frac{\nu \lambda_1}{4\big(1+\lambda_1\kappa\big)}}$ the following 
estimates hold true:
\begin{align}
  &\tau(t)\big(\|\nabla\bu_{t}\|^2+\kappa\,\|\tilde\Delta\bu_{t}\|^2\big)
+\nu\,e^{-2\alpha t}\,\displaystyle{\int_0^t}\bs(s)\,( \|\tilde\Delta\bu_s(s)\|^2+\|\nabla p_s(s)\|^2)ds\leq C,\label{a1}\\
  &e^{-2\alpha t}\displaystyle{\int_0^t}\bs(s)(\|\bu_{ss}(s)\|^2+\kappa\,\|\nabla\bu_{ss}(s)\|^2)ds\leq C,\label{a2}\\
  &e^{-2\alpha t}\displaystyle{\int_0^t}
  \bs(s)(\kappa\|\nabla\bu_{ss}(s)\|^2+\kappa^2\,\|\tilde\Delta\bu_{ss}(s)\|^2)ds\leq C.\label{a3r}
\end{align}
 Note that, 
here and everywhere else in the consecutive analysis the constant 
$C=C(\kappa,\nu,\alpha, \lambda_1, M)$ is independent of inverse powers of $\kappa$. 
\end{ldf}
\noindent{\it Proof.} Rewrite (\ref{2.3atw}) and differentiate the resulting equation with respect to time to arrive at
\begin{align}\label{a3}
 (\bu_{tt}, \bphi)-\kappa\,(\tilde\Delta\bu_{tt},\bphi)-\nu\, (\tilde\Delta\bu_t, \bphi)
 +(\bu_t\cdot\nabla \bu + \bu \cdot \nabla \bu_t, \bphi) =(\f_t, \bphi), \forall \bphi \in \bJ_1.
\end{align}
Choose $ \phi= -\bs(t)\tilde\Delta\bu_t$ in (\ref{a3}) to obtain
\begin{align}\label{2a}
 \frac{1}{2}\frac{d}{dt}\bs\left(\|\nabla\bu_{t}\|^2+\kappa\,\|\tilde\Delta\bu_{t}\|^2\right)
 &+\nu\,\bs\, \|\tilde\Delta\bu_t\|^2=C(\alpha,\lambda_1)e^{2\alpha t}\big(\|\nabla\bu_{t}\|^2
 +\kappa\,\|\tilde\Delta\bu_{t}\|^2\big)\nonumber\\
 &+\bs\,(\bu_t\cdot\nabla \bu,\tilde\Delta\bu_t) 
 + \bs\,(\bu \cdot \nabla \bu_t, \tilde\Delta\bu_t) -\bs\,(\f_t, \tilde\Delta\bu_t)\nonumber\\
 &=C(\alpha,\lambda_1)e^{2\alpha t}\big(\|\nabla\bu_{t}\|^2
 +\kappa\,\|\tilde\Delta\bu_{t}\|^2\big)+I_1+I_2+I_3, (\mbox{say}).
\end{align}
A use of Cauchy-Schwarz's inequality and Young's inequality lead to
\begin{align}\label{3a}
 |I_1|+|I_2|&\leq C
 \bs\,2^{1/2}\big(\|\bu_t\|^{1/2}\| \nabla \bu_t\|^{1/2} \|\nabla\bu\|^{1/2}\|\tilde\Delta\bu\|^{1/2}
 +\|\nabla\bu\|^{1/2}\|\tilde\Delta\bu\|^{1/2}\|\nabla\bu_t\|\big)\|\tilde\Delta\bu_t\|\nonumber\\
 &\leq C(\epsilon)\bs\,\big(\|\bu_t\|\| \nabla \bu_t\| \|\nabla\bu\|\|\tilde\Delta\bu\|
 +\|\nabla\bu\|\|\tilde\Delta\bu\|\|\nabla\bu_t\|^2\big)+\epsilon\,\bs\,\|\tilde\Delta\bu_t\|^2\nonumber\\
 &\leq  
 C(\epsilon,\lambda_1)\bs\,\|\nabla\bu_t\|^2\|\tilde\Delta\bu\|^2
 +\epsilon\,\bs\,\|\tilde\Delta\bu_t\|^2.
 \end{align}
Once again, apply Cauchy-Schwarz's inequality and Young's inequality to bound $|I_3|$ as
\begin{align}
 \label{4a}
|I_3|\leq C\,\bs\|\f_t\|\|\tilde\Delta\bu_t\|
\leq C(\epsilon)\bs\|\f_t\|^2+\epsilon\,\bs\,\|\tilde\Delta\bu_t\|^2.
\end{align}
After using (\ref{3a})-(\ref{4a}) in (\ref{2a}) with a proper choice of $\epsilon$, integrate the resulting equation with 
respect to time from $0$ to $t$ to arrive at
\begin{align}
 \label{6a}
 \bs\big(\|\nabla\bu_{t}\|^2&+\kappa\,\|\tilde\Delta\bu_{t}\|^2\big)
+\nu\,\displaystyle{\int_0^t}\bs(s)\, \|\tilde\Delta\bu_s(s)\|^2ds\leq C(\alpha,\lambda_1,\nu)\bigg(
\displaystyle{\int_0^t}e^{2\alpha s}\big(\|\nabla\bu_{s}(s)\|^2
 +\kappa\,\|\tilde\Delta\bu_{s}(s)\|^2\big)ds\nonumber\\
 &+\displaystyle{\int_0^t}\bs(s)\|\nabla\bu_s(s)\|^2\|\tilde\Delta\bu(s)\|^2\,ds+\displaystyle{\int_0^t}\,\bs(s)\|\f_s(s)\|^2
 ds\bigg).
\end{align}
Apply Lemma \ref{l1} and assumption ({\bf A2}) in (\ref{6a}). Then, multiply the resulting equation by $e^{-2\alpha t}$ 
to arrive at the desired {\it a priori} estimates of $\bu$ in (\ref{a1}).\\
\noindent
Next, differentiate (\ref{2.3atw}) and substitute $\bphi=\bs \bu_{tt}$ in the resulting equation to observe that
\begin{align}
 \label{8a1}
 \bs(\|\bu_{tt}\|^2+\kappa\,\|\nabla\bu_{tt}\|^2)=-\nu\,\bs\,(\nabla\bu_t, \nabla\bu_{tt})
 -\bs\,(\bu_t\cdot\nabla \bu + \bu \cdot \nabla \bu_t, \bu_{tt})+\bs\,(\f_t, \bu_{tt}).
\end{align}
After rewriting the first term on the right-hand side of (\ref{8a1}), apply Cauchy-Schwarz's inequality, Young's inequality and obtain
\begin{align}
 \label{9a}
 \bs(\|\bu_{tt}\|^2+\kappa\,\|\nabla\bu_{tt}\|^2)\leq 
 C(\nu)\bs\,(\|\tilde\Delta\bu_t\|^2+\|\nabla\bu_t\|^2+\|\tilde\Delta\bu\|^2+\|\f_t\|^2).
\end{align}
An integration of (\ref{9a}) with respect to time from $0$ to $t$, a multiplication by $e^{-2\alpha t}$ and a use of (\ref{a1}), 
Lemma \ref{l1}, assumption ({\bf A2}) complete the proof of (\ref{a2}).\\
\noindent
Now to derive (\ref{a3r}), substitute $\bphi=-\tilde \Delta \bu_{tt}$ in (\ref{2.3atw}) and use Cauchy-Schwarz's 
inequality, Young's inequality to yield
\begin{align}\label{18a}
 \|\nabla\bu_{tt}\|^2+\kappa\,\|\tilde\Delta\bu_{tt}\|^2\leq \frac{C(\nu)}{\kappa}
 \left( \|\tilde\Delta\bu_t\|^2+\|\nabla\bu_t\|^2+\|\tilde\Delta\bu\|^2+\|\f_t\|^2\right).
\end{align}
Multiply (\ref{18a}) by $\kappa\,\bs$ and integrate the resulting equation with respect to time from $0$ to $t$ to obtain
\begin{align}\label{20a}
 \displaystyle{\int_0^t}\bs(s)(\kappa\|\nabla\bu_{ss}(s)\|^2+\kappa^2\,\|\tilde\Delta\bu_{ss}(s)\|^2)ds\leq C
 \displaystyle{\int_0^t}\bs(s)\,\left(\nu\,\|\tilde\Delta\bu_s(s)\|^2+\|\nabla\bu_s(s)\|^2+\|\tilde\Delta\bu(s)\|^2+\|\f_s(s)\|^2\right)ds.
\end{align}
Multiply (\ref{20a}) by $e^{-2\alpha t}$ and use (\ref{a1}), Lemma \ref{l1} to arrive at the desired result in (\ref{a3r}).\\
Now to prove pressure estimate in (\ref{a1}), rewrite (\ref{2.2atw}). Then, differentiate the resulting 
equation with respect to time and obtain
\begin{align}\label{12a}
 ( \nabla p_t, \bphi)=(\bu_{tt}, \bphi) -\kappa (\tilde\Delta \bu_{tt}, \bphi) -\nu (\tilde\Delta \bu_t, \bphi)
 + ( \bu_t \cdot
\nabla \bu, \bphi)+( \bu \cdot
\nabla \bu_t, \bphi)-(\f_t,\bphi). 
\end{align}
Choose $\bphi=\nabla p_t$ in (\ref{12a}). Then, apply Cauchy-Schwarz's inequality and generalized H\"{o}lder's inequality to 
find that 
\begin{align}\label{13a}
 \|\nabla p_t\|\leq C(\|\bu_{tt}\|+\kappa\,\|\tilde\Delta\bu_{tt}\|+\nu\,\|\tilde\Delta\bu_t\|+\|\nabla\bu_t\|+\|\tilde\Delta\bu\|+\|\f_t\|).
\end{align}
After squaring both sides of (\ref{13a}), multiply it by $\bs$ and integrate with respect to time from $0$ to $t$ to arrive at
\begin{align}\label{14a}
 \displaystyle{\int_0^t}\bs(s)\,\|\nabla p_s(s)\|^2ds\leq C\displaystyle{\int_0^t}&\bs(s)
 (\|\bu_{ss}(s)\|^2+\kappa^2\,\|\tilde\Delta\bu_{ss}(s)\|^2+\nu\,\|\tilde\Delta\bu_s(s)\|^2\nonumber\\
 &+\|\nabla\bu_s(s)\|^2
 +\|\tilde\Delta\bu(s)\|^2+\|\f_s(s)\|^2)ds.
\end{align}
A use of estimates of $\bu$ from (\ref{a1}), (\ref{a2}), Lemma \ref{l1}, assumption ({\bf A2}) and a multiplication by $e^{-2\alpha t}$ lead to 
\begin{align}\label{15a}
e^{-2\alpha t} \displaystyle{\int_0^t}\bs(s)\,\|\nabla p_s(s)\|^2ds\leq C.
 \end{align}
 This completes the proof of Lemma \ref{l2}.\hfill{$\Box$}\\
 \noindent
To derive uniform estimates in time, we assume the following uniqueness condition: 
\begin{eqnarray}\label{tn}
 \frac{N}{\nu^2}\|\f\|_{L^{\infty}(0,\infty;\bL^2)}<1\,\,\, \rm{and}\,\,N=
 \sup_{\bu,\bv,\bw\in \bH_0^1(\Omega)}\frac{b(\bu,\bv,\bw)}{\|\nabla\bu\|\|\nabla\bv\|\|\nabla\bw\|}.
\end{eqnarray}
\section{\bf Semidiscrete Approximation}
\setcounter{tdf}{0}
\setcounter{ldf}{0}
\setcounter{cdf}{0}
\setcounter{equation}{0}
Let ${\bf H}_h$ and $L_h$ be the finite-dimensional subspaces of 
${\bf H}_0^1 $ and $L^2$, respectively, such that, there exist operators $i_h$ and $j_h$ satisfying the following 
approximation properties:  \\
\noindent
({\bf B1}). For each $\bw \in {\bf {J}}_1 \cap {\bf {H}}^2 $ and $ q \in H^1/ {\rm I\!R}$, there exist approximations 
$i_h \bw \in {\bf
{J}}_h $ and $ j_h q \in L_h $ such that
\be
\|\bw-i_h\bw\|+ h \| \nabla (\bw-i_h \bw)\| \le K_0 h^2
\| \bw\|_2, \;\;\;\; \| q - j_h q
 \|_{L^2 /{\rm I\!R}} \le K_0 h \| q\|_{H^1 / {\rm I\!R}}.
\ee
\noindent
Note that, $h>0$ be a discretization parameter with $0<h<1$.\\
\noindent
Here, it can be noted that the operator $b(\cdot, \cdot, \cdot)$ preserves the
antisymmetric properties of the original nonlinear term, i.e.,
\begin{align}\label{eq1}
b(\bv_h, \bw_h, \bw_h) = 0 \;\;\; \forall \bv_h, \bw_h \in
{\bH}_h.
\end{align}
\noindent
The discrete analogue of the weak formulation (\ref{2.2atw}) is as follows:\\
\noindent
 Find $\bu_h(t) \in {\bf H}_h$ and $p_h(t) \in L_h$ such that $ \bu_h(0)
 = \bu_{0h} $ and for $t>0$,
\begin{align}
\label{3.1}
(\bu_{ht}, {\bphi}_h)+\kappa a ( \bu_{ht}, {\bphi}_h) +\nu a ( \bu_h, {\bphi}_h) &+ b( \bu_h, \bu_h, {\bphi}_h)
-(p_h, \nabla \cdot {\bphi}_h)
=0 \;\;\; \forall {\bphi}_h \in {\bf H}_h,\nonumber\\
(\nabla \cdot \bu_h, \chi_h) &= 0 \;\;\; \forall \chi_h \in L_h,
\end{align}
\noindent
where $\bu_{0h} \in {\bf H}_h $ is a suitable approximation of $\bu_0\in
{\bf J}_1$.\\
\noindent
For subsequent analysis, we define a suitable approximation of $\bJ_1$ by introducing the discrete 
incompressibility condition in $\bH_h$ and
call the resulting subspace as $\bJ_h$. Thus, $\bJ_h$ is defined as
\be
{\bf J}_h = \{ \bv_h \in {\bf H}_h : (\chi_h, \nabla \cdot \bv_h)
=0 \;\;\; \forall \chi_h \in L_h \}.
\ee
 Note that, the space $\bJ_h$ is not a  subspace of $\bJ_1$. 
Now, an equivalent form of (\ref{3.1}) is defined as:\\
\noindent
Find $\bu_h (t) \in {\bf J}_h $ such that $\bu_h(0) = \bu_{0h} $
and for $t>0$,
\ben
\label{3.2}
(\bu_{ht},{\bphi}_h) +\kappa a (\bu_{ht},{\bphi}_h)
+\nu a (\bu_h,{\bphi}_{h})
 + b( \bu_h, \bu_h, {\bphi}_h)=0 
 \;\;\; \forall {\bphi}_h \in {\bf J}_h.
\een
 For proof of the global existence of a unique solution of (\ref{3.2}), one may refer to \cite{bnpdy}.\\
\noindent
In order to deal with the pressure estimates in subsequent analysis, 
we assume the pair $({\bf H}_h,{L_h/N_h})$ satisfies a uniform inf-sup condition:\\
\noindent
({\bf B2}). For every $q_h \in L_h$, there exist a non-trivial
function $\bphi_h \in \bf H_h$ and a positive constant $K_1$, independent of $h$, such that, 
\be
|(q_h, \nab \cdot {\bphi}_h)| \ge K_1\|\nab {\bphi}_h \| \| q_h\|_{
L^2/N_h}.
\ee
 The following properties of the $L^2$ projection $P_h:\bL^2\rightarrow \bJ_h$ can be derived using conditions ({\bf B1})-({\bf B2}) (
 for a proof, see (\cite{GR}, \cite{HR82}):  
\ben
\label{3.4}
\|\bphi- P_h \bphi\|+ h \|\nabla P_h \bphi\| \leq C h
\|\nabla \bphi\|\,\,\,\,\,\,\,\,\forall {\bphi \in \bJ_1},
\een
and
\ben
\label{3.5}
\|\bphi- P_h \bphi\|+ h \|\nabla (\bphi-P_h \bphi)\| \leq C h^2
\|\tilde \Delta \bphi\|\,\,\,\,\forall\bphi \in \bJ_1 \cap \bH^2.
\een
We may define the discrete  operator $\Delta_h: \bH_h
\rightarrow \bH_h$ through the bilinear form $a (\cdot, \cdot)$ as
\begin{eqnarray}
\label{3.6}
a(\bv_h, {\bphi}_h) = (- \Delta_h \bv_h, \bphi_h) \;\;\;\;
\forall \bv_h, {\bphi}_h \in \bH_h.
\end{eqnarray}
Set the discrete analogue of the Stokes operator $\tl{\Delta} =
P \Delta $ as
 $\tl{\Delta}_h = P_h \Delta_h $.
Examples of subspaces $\bf H_h$ and $L_h$ satisfying assumptions $(\rm{\bf {B1}})$ and $(\rm{\bf {B2}})$ 
in the context of both conforming and non-conforming analysis can be 
found in \cite{BF}, \cite{BP} and \cite{HR82}.\\
\noindent
We recall below in Lemma \ref{L72}, some {\it a priori} bounds of $\bu_h$ which will
be used in the derivation of fully discrete error estimates in the subsequent section. For proof, one may refer to 
\cite{ANS6} (Lemma 4.2), \cite{AKP} (Lemma 3.2).
\begin{ldf} \label{L72}
Let the assumptions {\rm{({\bf {A1}})}}-{\rm{({\bf {A2}})}} hold. Then, there exists a 
positive constant $C =C(\kappa,\nu,\alpha, \lambda_1, M)$ such that 
for $\displaystyle{0\leq\alpha< \frac{\nu \lambda_1}{4\big(1+\lambda_1\kappa\big)}}$ the following 
estimates hold true: 
\begin{align}
&\|\bu_{ht}(t)\|^2+\|\tl{\Delta}_h\bu_h(t)\|^2+
e^{-2\alpha t}\int_0^t e^{2\alpha s}(\|\nabla\bu_h(s)\|^2+\|\tl{\Delta}_h\bu_h(s)\|^2+\|\nab\bu_{hs}(s)\|^2)ds \le C,
\nonumber\\
&e^{-2\alpha t}\int_0^t e^{2\alpha s}(\|\bu_{hss}(s)\|_{-1}^2+\kappa\,\|\nabla\bu_{hss}(s)\|^2)\,ds \le C.\nonumber
\end{align}
\end{ldf}
\hfill {$\Box$}\\ 
\noindent
Now, in Theorem \ref{T31*}, the main results of the article are stated in which we present the semidiscrete optimal error 
estimates of the velocity and pressure. The proofs are established in Sections ${\bf 5}$.
\begin{tdf}\label{T31*}
Let the assumptions {\rm{({\bf {A1}})}}-{\rm{({\bf {A2}})}} 
and {\rm (\bf{B1})-(\bf{B2})} 
be satisfied. Let $\bu_{0h}=P_h \bu_0$, then, there exists a positive constant $C$ depending on $\kappa$, $\lambda_1$, 
$\nu$, $\alpha$ and $M$, 
such that, for fixed $T > 0$ with $t\in(0,\,T)$ and for 
$\displaystyle{0\leq\alpha < \frac{\nu \lambda_1}{4\big(1+\lambda_1\kappa\big)}}$,
the following estimates hold true:
\begin{align*}
\|(\bu-\bu_h)(t)\|+h\|\nab (\bu-\bu_h)(t)\|\le K(t) h^2,
\end{align*}
$K(t)=C e^{Ct}$. Under the uniqueness condition (\ref{tn}), $K(t)=C$, that is, the estimates are uniform in time.
\end{tdf}
\begin{tdf}\label{T32*}
Under the hypotheses of Theorem \ref{T31*},  
there exists a positive constant $C$ depending on $\kappa$, $\lambda_1$, $\nu$, $\alpha$ and $M$, 
such that, for all $t>0$, the following holds true:
$$
\|(p-p_h)(t)\|_{L^2/N_h}\leq K(t) h.
$$
\end{tdf}
\noindent
Here again, under the uniqueness condition (\ref{tn}), $K(t)=C$, that is, the estimate holds uniformly with respect to time. 
\section{\bf Semidiscrete Finite Element Error Estimates}
\setcounter{tdf}{0}
\setcounter{ldf}{0}
\setcounter{cdf}{0}
\setcounter{equation}{0}
This section deals with the optimal error estimates of velocity and pressure. Note that, 
since $\bJ_h$ is not a subspace of $\bJ_1$, the weak solution $\bu$ satisfies
\begin{eqnarray}\label{E62}
(\bu_t,{\bphi}_h)+\kappa a(\bu_{t},{\bphi}_h)+\nu a(\bu,{\bphi}_h)=-b(\bu,\bu,{\bphi}_h)+(p,\nabla \cdot {\bphi}_h)
+(\f,\bphi_h) \,\,\,\, \forall {\bphi}_h \in \bJ_h.  
\end{eqnarray}
\noindent
Set ${\bf e}=\bu-\bu_h$. Then, subtract (\ref{E62}) from (\ref{3.2}) to arrive at
\begin{eqnarray}\label{E63}
({\bf e}_{t},{\bphi}_h)+\kappa a({\bf e}_{t},{\bphi}_h)+\nu a({\bf e},{\bphi}_h)= {\bf \Lambda}({\bphi}_h)
+(p,\nabla \cdot {\bphi}_h),
\end{eqnarray}
\noindent
where ${\bf \Lambda}({\bphi}_h)=-b(\bu,\bu,{\bphi}_h)+b(\bu_h,\bu_h,{\bphi}_h)$.
Below, we derive the optimal error estimates of $||\e(t)||$ and $||\nabla\e(t)||$, for $t>0$.\\
\noindent
In order to deal with the nonlinearity, an intermediate solution $\bv_h$ is introduced which is a finite element Galerkin 
approximation 
to a linearized Kelvin-Voigt equation. The solution $\bv_h$ satisfies 
\ben
\label{E614}
 (\bv_{ht}, {\bphi}_h)+\kappa\,a(\bv_{ht},{\bphi}_h)+\nu a (\bv_h, {\bphi}_h)
=-b(\bu,\bu,{\bphi}_h)
\;\;\;\; \forall \bphi \in {\bf J}_h,
\een
with $\bv_h(0)=P_h \bu_0.$\\
\noindent
Now, we split ${\bf {e}} $ as
\be
{\bf  {e}} := \bu - \bu_h = (\bu - \bv_h) + ( \bv_h - \bu_h)
= \bxi + \bta.
\ee
Here, $\bxi$ is the error due to the approximation using a linearized Kelvin-Voigt equation (\ref{E614}), whereas 
$\bta$ denotes the error
due to the non-linearity in the equation.
A subtraction of (\ref{E614}) from (\ref{E62}) leads to the equation in $\bxi$ as
\begin{eqnarray}\label{E615}
(\bxi_t,{\bphi}_h)+\kappa\,a({\bxi}_t,{\bphi}_h)+\nu a(\bxi,{\bphi}_h)=(p,\nabla \cdot {\bphi}_h)\;\;\;\;\forall {\bphi}_h\in {\bf {J}}_h.
\end{eqnarray}
In order to derive optimal error estimates of $\bxi$ in $L^{\infty}(\bL^2)$ and
$L^{\infty}(\bH^1)$-norms, we introduce the following
auxiliary  projection $V_h$ such that $V_h\bu:[0,\infty)\rightarrow \bJ_h$
satisfying
\begin{eqnarray}\label{E616}
\kappa  a(\bu_t-V_h\bu_t,{\bphi}_h)+ \nu a(\bu-V_h\bu,{\bphi}_h)=(p,\nabla \cdot {\bphi}_h)\; \; \;
\forall {\bphi}_h\in \bJ_h,
\end{eqnarray}
where $V_h\bu(0)=P_h \bu_0.$\\
\noindent
With $V_h\bu $ defined as above, we now split $\bxi$ as
$$
\bxi:=(\bu-V_h\bu)+(V_h\bu-\bv_h)=\bzeta+\brho.
$$
To obtain estimates for ${\e}$, first of all, we establish a few estimates of $ \bzeta $ 
in Lemmas \ref{L61}-\ref{La4}. Then with the help of $\bzeta$ estimates,  we derive various estimates of $\brho$ 
and $\nabla \brho$ 
 in Lemmas \ref{L65*} and \ref{L63}.
 Finally, in Lemma \ref{L63n}, we derive estimates for $\bta$ and complete the proof of Theorem \ref{T31*}.
\begin{ldf}\label{L61}
Assume that assumptions {\rm{({\bf {A1}})}}-{\rm{({\bf {A2}})}}  and {\rm{({\bf B1})-({\bf B2})}} are satisfied. Then, 
there exists a positive 
constant $C=C(\kappa,\nu,\alpha,\lambda_1,M)$ such that for $\displaystyle{0\leq\alpha< \frac{\nu\lambda_1}{4\big(1+\kappa\lambda_1\big )}}$, 
the following estimate holds true:
\begin{eqnarray*}
\bs^{-1}(t)\int_0^t e^{2\alpha s}\|\nabla(\bu-V_h\bu)(s)\|^2 ds\leq C  h^2.
\end{eqnarray*}
\end{ldf}
\noindent{\it Proof.} 
 Multiply (\ref{E616}) by $e^{\alpha t}$ with $\bzeta =\bu-V_h\bu$, use $e^{\alpha t}\bzeta_t= 
 \hat\bzeta_t-\alpha \hat\bzeta$ and 
 substitute ${\bphi}_h=P_h\hat{\bzeta}=\hat\bzeta+(P_h \hat\bu-\hat \bu)$ to arrive at   
\begin{align}\label{E619}
 \kappa \frac{d}{dt}\|\nabla \hat{\bzeta}\|^2 +2(\nu-&\kappa \alpha)\|\nabla \hat{\bzeta}\|^2 = 2\kappa \frac{d}{dt}
a ( \hat\bzeta, \hat\bu-P_h \hat{\bu})-2\kappa\,a (\hat\bzeta, \frac{d}{dt}(\hat\bu-P_h \hat{\bu}))\nonumber\\
&+ 2(\nu-\kappa \alpha)\, a(\hat{\bzeta},
\hat{\bu}-P_h\hat{\bu})
+ 2(\hat{p}-j_h\hat{p},
\nabla \cdot P_h\hat{\bzeta}).
\end{align}
Integrate (\ref{E619}) with respect to time from $0$ to $t$ and 
apply (\ref{3.4}) along with Young's inequality. A simplification of resulting equation with a use of 
$\|\bu_0-P_h\bu_0\|=\|\bu(0)-P_h\bu(0)\|$ yields 
\begin{align}
\label{Ea3}
\kappa\|\nabla \hat{\bzeta}\|^2& +2(\nu-\kappa \alpha)
\displaystyle{\int_0^t}\|\nabla \hat{\bzeta}(s)\|^2ds \leq 2\kappa a(\hat\bzeta,\hat \bu-P_h\hat \bu)-
\kappa \|\nabla(\bu(0)-P_h\bu(0))\|^2\nonumber\\
&+2\kappa\,\displaystyle{\int_0^t}\,
\|\nabla\hat\bzeta(s)\|\|\nabla(\hat{\bu}_s-P_h \hat{\bu}_s)(s)\|ds+ 2(\nu-\kappa \alpha)\,\displaystyle{\int_0^t} 
\|\nabla\hat{\bzeta}(s)\|\|\nabla(\hat{\bu}-P_h\hat{\bu})(s)\|ds\nonumber\\
&+ 2\,\displaystyle{\int_0^t}\,\|(\hat{p}-j_h\hat{p})(s)\|
\|\nabla \cdot P_h\hat{\bzeta}(s)\|ds.
\end{align}
After applying Cauchy-Schwarz's inequality in the first term of right-hand side, 
use Young's inequality with $p=1/2$, $q=1/2$ to obtain
\begin{align}\label{Ea1}
2\kappa a(\hat\bzeta,\hat \bu-P_h\hat \bu)&\leq \kappa\|\nabla\hat\bzeta\|^2+\kappa\|\nabla(\hat \bu-P_h\hat \bu)\|^2.
\end{align}
A use of (\ref{Ea1}) in (\ref{Ea3} ) leads to
\begin{align}
\label{Ea2}
\kappa\|\nabla \hat{\bzeta}\|^2& +2(\nu-\kappa \alpha)\displaystyle{\int_0^t}\|\nabla \hat{\bzeta}(s)\|^2ds 
\leq \kappa\|\nabla\hat\bzeta\|^2+(\kappa\|\nabla(\hat \bu-P_h\hat \bu)\|^2\nonumber\\
&-\kappa\|\nabla(\bu(0)-P_h\bu(0))\|^2)+2\kappa\,\displaystyle{\int_0^t}\,
\|\nabla\hat\bzeta(s)\|\|\nabla(\hat{\bu}_s-P_h \hat{\bu}_s)(s)\|ds\nonumber\\
&+ 2(\nu-\kappa \alpha)\,\displaystyle{\int_0^t} \|\nabla\hat{\bzeta}(s)\|\|\nabla(\hat{\bu}-P_h\hat{\bu})(s)\|ds
+ 2\,\displaystyle{\int_0^t}\,\|(\hat{p}-j_h\hat{p})(s)\|
\|\nabla \cdot P_h\hat{\bzeta}(s)\|ds.
\end{align}
The first term on both sides will cancel out. To deal with the second term on right-hand side, rewrite it as
\begin{align}
\label{Ea4}
\kappa&\|\nabla(\hat \bu-P_h\hat \bu)\|^2-\kappa\|\nabla(\bu(0)-P_h\bu(0))\|^2=
\kappa \displaystyle{\int_0^t}\frac{d}{ds}\left\|\nabla(\hat{\bu}-P_h\hat{\bu})\right\|^2ds\nonumber\\
&=\kappa  \displaystyle{\int_0^t}\frac{d}{ds}e^{2\alpha s}(\nabla(\bu-P_h\bu),\nabla(\bu-P_h\bu))ds\nonumber\\
&\leq C(\alpha) h^2\displaystyle{\int_0^t}e^{2\alpha s}
\left(\kappa\|\tilde{\Delta}\bu_s\|\|\tilde{\Delta}\bu\|+\kappa\|\tilde{\Delta}\bu\|^2\right)ds.
\end{align}
Apply (\ref{Ea4}) in (\ref{Ea2}) along	with Young's inequality, (\ref{3.5}) and {\rm{({\bf B1})}} to arrive at
\begin{align}
(\nu-\kappa \alpha)&\int_0^t\|\nabla \hat\bzeta\|^2 ds \leq 
C(\kappa,\nu,\alpha) h^2\int_0^t e^{2\alpha s}( \kappa\|\tilde{\Delta}\bu_s\| ^2 +\kappa\|\tilde{\Delta}\bu\|^2
+ \nu\|\tilde{\Delta} \bu \|^2+ \|\nabla p\|^2)ds.\nonumber
\end{align}
A use of {\it a priori} bounds for $\bu$ and $p$ stated in Lemma \ref{l1} completes the proof.
\hfill{$\Box$}
\begin{remark}\label{lr1}
 Note that using (\ref{Ea4}), we rewrite the second term on the right-hand side of (\ref{Ea2}) and thereby 
 write the entire right-hand side of (\ref{Ea2}) as an integration. This plays an important role 
 in achieving weight function $\bs$ in the desired estimates. The presence of $\bs$ in the estimates is crucial 
 in order to deal with the regularity issues at $t=0$ while making error estimates independent of $\kappa$. 
    \end{remark}
\noindent
Next, we prove the estimates for the time derivative of $\bzeta$.
\begin{ldf}\label{L62}
Under the assumptions {\rm{({\bf {A1}})}}-{\rm{({\bf {A2}})}}  and {\rm{({\bf B1})-({\bf B2})}}, 
there exists a positive constant $C=C(\kappa,\nu,\alpha,\lambda_1,M)$ such that for 
$\displaystyle{0\leq\alpha< \frac{\nu\lambda_1}{4\big(1+\kappa\lambda_1\big)}}$, the following estimates hold true:
\begin{eqnarray}
 &&\bs^{-1}(t)\kappa\int_0^t\bs\|\nabla (\bu_s(s)-V_h\bu_s(s))\|^2 ds
 +\nu \|\nabla (\bu-V_h\bu)(t)\|^2 \leq C h^2,\label{ee10}\\
 &&\bs^{-1}(t)\kappa^2\displaystyle{\int_0^t} e^{2\alpha s}\|\nabla(\bu_s(s)-V_h\bu_s(s))\|^2ds\leq C h^2.\label{ee2} 
\end{eqnarray}
\end{ldf}
\noindent{\it Proof}. Recall (\ref{E616}) now with ${\bphi}_h= P_h{\bzeta}_t=\bzeta_t+(P_h \bu_t-\bu_t)$ to find that
\begin{align}\label{E640}
 2\kappa\|\nabla \bzeta_t\|^2+\nu\frac{d}{dt}\|\nabla {\bzeta}\|^2=2( p, \nabla \cdot  P_h\bzeta_t)
 +2\kappa\,a(\bzeta_t,(\bu_t-P_h \bu_t))+2\nu a( \bzeta, (\bu_t-P_h \bu_t)).
\end{align}
Rewrite the first term on the right-hand side of (\ref{E640}) as 
\begin{align}
2( p, \nabla \cdot  P_h\bzeta_t)=2\frac{d}{dt}(p, \nabla \cdot  P_h\bzeta)-2(p_t,\nabla \cdot  P_h\bzeta)
\end{align}
and substitute in (\ref{E640}) to obtain
\begin{align}
\label{Ea5}
2\kappa\|\nabla \bzeta_t\|^2&+\nu\frac{d}{dt}\|\nabla {\bzeta}\|^2=2\frac{d}{dt}(p, \nabla \cdot  P_h\bzeta)-2(p_t,\nabla \cdot  P_h\bzeta)\nonumber\\
&+2\kappa\,a(\bzeta_t,(\bu_t-P_h \bu_t))+2\nu a( \bzeta, (\bu_t-P_h \bu_t)).
\end{align}
After using the Cauchy-Schwarz inequality and discrete incompressibility condition  in (\ref{Ea5}), multiply the 
resulting equation by $\bs(t)$  
to arrive at 
\begin{align}\label{Ea6}
2\kappa &\bs\|\nabla \bzeta_t\|^2+\nu\frac{d}{dt}\bs\|\nabla{\bzeta}\|^2\leq \nu \bs_t\|\nabla \bzeta\|^2 +2\frac{d}{dt}
\bs (p-j_h p, \nabla \cdot  P_h\bzeta)-2\bs_t (p-j_h p, \nabla \cdot  P_h\bzeta)\nonumber\\
&+2\bs\|p_t-j_h p_t\|\|\nabla \cdot  P_h\bzeta\|
+2\kappa\bs \| \nabla \bzeta_t\|\|\nabla(\bu_t-P_h \bu_t)\|
+ 2\nu \bs\|\nabla \bzeta\|\| \nabla(\bu_t-P_h \bu_t)\|.
\end{align}
Integrate (\ref{Ea6}) with respect to time from $0$ to $t$ and apply Young's inequality, (\ref{3.5}), $\rm({\bf B1})$ 
to obtain
\begin{align}
\label{Ea8}
\kappa \displaystyle{\int_0^t}\bs&\|\nabla \bzeta_s\|^2ds+\nu\bs\|\nabla{\bzeta}\|^2\leq C(\nu)\bigg( \displaystyle{\int_0^t}\bs_{s}\|\nabla \bzeta\|^2ds +\bs h^2\|\nabla p\|^2 \nonumber\\
&+h^2\displaystyle{\int_0^t}\bs_s \|\nabla p\|^2ds+\displaystyle{\int_0^t}\bs_s\| \nabla \bzeta\|^2ds
+h^2\displaystyle{\int_0^t}\frac{\bs^2}{\bs_s}\|\nabla p_s\|^2ds\nonumber\\
&+\kappa h^2\displaystyle{\int_0^t}\bs\|\tilde\Delta\bu_s\|^2ds
+ \nu\displaystyle{\int_0^t} \frac{\bs^2}{\bs_s}\|\tilde\Delta\bu_s\|^2ds\bigg).
\end{align}
Now, the desired results in (\ref{ee10}) follows by using Lemmas \ref{l1} and \ref{L61} in (\ref{Ea8}). \\
\noindent
Next to prove (\ref{ee2}), substitute $\bphi_h=P_h\bzeta_t$ with $\bzeta=\bu-V_h\bu$ in (\ref{E616}) and arrive at
\begin{align}
\label{ea11}
\kappa\,\|\nabla\bzeta_t\|^2=\kappa\,a(\bzeta_t,\bu_t-P_h\bu_t)-\nu\,a(\bzeta,P_h\bzeta_t)+(p,\nabla\cdot P_h\bzeta_t).
\end{align}
A use of Cauchy-Schwarz's inequality, (\ref{3.5}) and discrete incompressibility condition in (\ref{ea11}) yield
\begin{align}
\label{ea13}
\kappa\,\|\nabla\bzeta_t\|\leq C(h\kappa\,\|\tilde{\Delta}\bu_t\|+\nu\,\|\nabla \bzeta\|+h\|\nabla p \|).
\end{align}
After taking a square of (\ref{ea13}) on both sides, multiply the resulting equation by $e^{2\alpha t}$. Then, integrate 
with respect to time 
from $0$ to $t$ and use bounds from Lemmas \ref{l1}, \ref{L61} to arrive at the desired result.
This completes the rest of the proof.
\hfill{$\Box$}\\
\begin{ldf}\label{La3}
  	Under the assumptions {\rm{({\bf {A1}})}}-{\rm{({\bf {A2}})}}  and {\rm{({\bf B1})-({\bf B2})}}, there exists a positive constant 
  	$C=C(\kappa,\nu,\alpha,\lambda_1,M)$ such that for $\displaystyle{0\leq\alpha< \frac{\nu}{4\big(1+\kappa\lambda_1\big)}}$, 
  	the following estimate holds true:
  	\begin{eqnarray}
  \tau(t)\kappa\|\nabla(\bu_t-V_h\bu_t)(t)\|^2+\nu\bs^{-1}(t)\displaystyle{\int_0^t}\bs_1(s) \|\nabla(\bu_s-V_h\bu_s)(s)\|^2ds
  \leq C h^2.\nonumber
  	\end{eqnarray}
Here, $\tau(t):=\min\{t,\,1\}$ and $\bs_1(t):=\tau^2(t)\,e^{2\alpha t}$.
  \end{ldf}
  \noindent{\it Proof.} Differentiate (\ref{E616}) with respect to time and substitute $\bphi_h=P_h\bzeta_t$ to observe
  \begin{align}\label{ea41*}
  \frac{\kappa}{2}\,\frac{d}{dt}\|\nabla\bzeta_{t}\|^2+\nu\,\|\nabla\bzeta_t\|^2
  =\nu\,a(\bzeta_t,\bu_t-P_h\bu_t)+
   \kappa\,a(\bzeta_{tt},\bu_t-P_h\bu_t)+(p_t,\nabla\cdot P_h\bzeta_t).
  \end{align}
  Note that,  a use of
  \begin{align}\label{ea42}
  \kappa\,\|\nabla\bzeta_{tt}\|\leq C(h\kappa\,\|\tilde{\Delta}\bu_{tt}\|+\nu\,\|\nabla \bzeta_t\|+h\|\nabla p_t \|),
  \end{align}
 {\rm{({\bf B1}), (\ref{3.5}), Cauchy-Schwarz's inequality and Young's inequality in (\ref{ea41*}) lead to
  \begin{align}\label{ea43}
  \kappa\,\frac{d}{dt}\|\nabla\bzeta_{t}\|^2+\nu\,\|\nabla\bzeta_t\|^2\leq C(\nu)h(\|\tilde{\Delta}\bu_t\|^2
  + \kappa\,\|\tilde{\Delta}\bu_{tt}\|^2+\|\nabla p_t\|^2).
  \end{align}
  Multiplication of (\ref{ea43}) by $\bs_1$ and integration of the resulting equation from $0$ to $t$ yield
  \begin{align}
 \kappa\,\bs_1\|\nabla\bzeta_{t}\|^2+\nu\,\displaystyle{\int_0^t}\bs_1\,\|\nabla\bzeta_s\|^2ds\leq \kappa\displaystyle{\int_0^t}\bs_{1,s}
 \|\nabla\bzeta_{s}\|^2ds+ C(\nu)h^2 \displaystyle{\int_0^t}\bs_1 (\|\tilde{\Delta}\bu_s\|^2
 + \kappa\,\|\tilde{\Delta}\bu_{ss}\|^2+\|\nabla p_s\|^2)ds.\nonumber
  \end{align}
  A use of estimates from Lemmas \ref{l2}, \ref{L62} and a multiplication of resulting equation by $\bs^{-1}(t)$ complete 
  the proof of Lemma \ref{La3}.\hfill{$\Box$}\\
\noindent
Below, in Lemma \ref{L62*} we discuss the $L^2(\bL^2)$-estimate of $\bzeta$. The similar kind of estimate 
has already been 
discussed in Lemma 5.3 of \cite{ANS6}. The difference between the estimate of $\bzeta$ in Lemma 5.3 of \cite{ANS6} and 
Lemma \ref{L62*} in this article is the presence of weight function 
$\bs$ in Lemma \ref{L62*} which will be very helpful in making the error estimates independent of $\kappa$. 
Therefore, in order to justify the presence of $\bs$, we present a short proof highlighting only the modifications.
\begin{ldf}\label{L62*}
Under the assumptions {\rm{({\bf {A1}})}}-{\rm{({\bf {A2}})}}  and {\rm{({\bf B1})-({\bf B2})}}, there exists a positive 
constant $C=C(\kappa,\nu,\alpha,\lambda_1, M)$ such that 
for $\displaystyle{0\leq\alpha< \frac{\nu\lambda_1}{4\big(1+\kappa\lambda_1\big)}}$, 
the following estimate holds true for $t >0$:
\begin{eqnarray}
\bs^{-1}(t)\int_0^t e^{2\alpha s}\|(\bu-V_h\bu)(s)\|^2 ds \leq C h^4.\nonumber
\end{eqnarray}
\end{ldf}
\noindent{\it Proof.} For obtaining the desired estimates of $\bzeta$, we appeal to the Aubin-Nitsche duality
argument by assuming $(\bw,q)$ to be the unique
solution of the steady state Stokes system:
\begin{align}
-\nu\Delta \bw+\nabla q&=\hat{\bzeta}\ \ {\rm in}\ \ \Omega, \label{E624}\\
\nabla \cdot \bw&=0\ \ \ \ {\rm in}\ \ \Omega,  \label{E625} \\
\bw|_{\partial\Omega}&=0  \label{E626}
\end{align}
satisfying the following regularity result:
\begin{eqnarray}\label{E6281*}
\|\bw\|_2+\|q\|_{H^1/{\rm I\!R}}\le C\|\hat{\bzeta}\|.
\end{eqnarray}
Form $L^2$-inner product between (\ref{E624}) and $\hat{\bzeta}$ and use discrete incompressibility condition. Then, 
apply (\ref{E616}) with ${\bphi}_h$ replaced by $P_h\bw$ to obtain
\begin{align}\label{E62810*}
\|\hat{\bzeta}\|^2&= \nu\, a (\bw-P_h \bw, \hat \bzeta)-(q-j_{h}q,\nabla\cdot \hat \bzeta)
+(\hat p-j_h \hat p,\nabla \cdot (P_h \bw-\bw))\nonumber\\
&\qquad-\kappa\,a(e^{\alpha t} \bzeta_t, P_h \bw-\bw)-
   \kappa\,a(e^{\alpha t} \bzeta_t,\bw)
\end{align}
 Once again, form an $L^2$-inner product between (\ref{E624}) and $e^{\alpha t}{\bzeta}_t$ and use it to replace the 
 last term in (\ref{E62810*}) as follows
 \begin{align}\label{E6281**}
\|\hat{\bzeta}\|^2+\frac{\kappa}{2\nu}\frac{d}{dt}\|\hat\bzeta\|^2&= \frac{\alpha \kappa}{\nu}\|\hat\bzeta\|^2+\nu a (\bw-P_h \bw, \hat \bzeta)-(q-j_{h}q,\nabla\cdot \hat \bzeta)+(\hat p-j_h \hat p,\nabla \cdot (P_h \bw-\bw))\nonumber\\
&-\kappa\,a(e^{\alpha t}\bzeta_t, P_h \bw-\bw)-\frac{\kappa}{\nu}(q-j_h q,  e^{\alpha t}\nabla\cdot \bzeta_t).
\end{align}
Apply Cauchy-Schwarz's inequality, assumption ({\bf B1}) and regularity estimates (\ref{E6281*}) along with Young's inequality 
in (\ref{E6281**}). Then, integrate the resulting equation with respect to time from $0$ to $t$ to obtain 
\begin{align}
\label{ea19*}
&\frac{(\nu-\alpha \kappa)}{2\nu}\displaystyle{\int_0^t} \|\hat{\bzeta}\|^2ds+\frac{\kappa}{2\nu}\|\hat\bzeta\|^2\leq 
\frac{\kappa}{2\nu}\|\bzeta(0)\|^2+C(\kappa,\nu,\alpha)\,h^2 \displaystyle{\int_0^t}
\big(\| \nabla\hat \bzeta\|^2
+h^2\|\nabla\hat p\|^2+\kappa h^2\,\|e^{\alpha s}\tilde{\Delta} \bu_s\|^2\big)ds.
\end{align}
Using $\|\bzeta(0)\|^2=\|\bu_0-P_h\bu_0\|^2$, we write 
\begin{align}
\label{ea14}
\|\bzeta(0)\|^2&=\|\bu_0-P_h\bu_0\|^2= \|\hat\bu-P_h\hat\bu\|^2-\displaystyle{\int_0^t}
\frac{d}{ds}\left(e^{2\alpha s}\|\bu-P_h\bu\|^2\right)ds\nonumber\\
&=\|\hat\bu-P_h\hat\bu\|^2-\displaystyle{\int_0^t}e^{2\alpha s}\left((\bu-P_h\bu,\bu_s-P_h\bu_s)+(\bu_s-P_h\bu_s,\bu-P_h\bu)\right)ds\nonumber\\
&+2\alpha\displaystyle{\int_0^t}e^{2\alpha s}\|\bu-P_h\bu\|^2ds.
\end{align}
A use of orthogonality property of $P_h$ and Cauchy-Schwarz's inequality yield 
\begin{align}
\label{ea15}
\|\hat{\bu}-P_h\hat\bu\|^2&=(\hat\bu-P_h\hat\bu,\hat\bu-P_h\hat\bu)\nonumber\\
&=(\hat\bu-V_h\hat\bu,\hat\bu-P_h\hat\bu)\leq \|\hat\bu-V_h\hat\bu\|\|\hat\bu-P_h\hat\bu\|.
\end{align}
A simplification of (\ref{ea15}) leads to
\begin{align}
\label{ea16}
\|\hat{\bu}-P_h\hat\bu\|\leq \|\hat\bu-V_h\hat\bu\|=\|\hat\bzeta\|.
\end{align}
An application of (\ref{ea14}) and (\ref{ea16}) in (\ref{ea19*}) yield
\begin{align}
\label{ea20}
\frac{(\nu-\alpha \kappa)}{\nu}\int_0^t \|\hat{\bzeta}\|^2ds&\leq 
-\frac{\kappa}{2\nu}\displaystyle{\int_0^t}
\bigg(e^{2\alpha s}(\bu-P_h\bu,\bu_s-P_h\bu_s)\nonumber\\
&+e^{2\alpha s}(\bu_s-P_h\bu_s,\bu-P_h\bu)\bigg)ds+\frac{\alpha\kappa}{\nu}\displaystyle{\int_0^t}e^{2\alpha s}\|\bu-P_h\bu\|^2ds\nonumber\\
& +C(\kappa,\nu,\alpha)\,h^2 \displaystyle{\int_0^t}
\big(\| \nabla\hat \bzeta\|^2
+h^2\|\nabla\hat p\|^2+\kappa h^2\,\|e^{\alpha s}\tilde{\Delta} \bu_s\|^2\big)ds.
\end{align}
By using (\ref{3.5}), we arrive at
\begin{align}
\label{E6282}
(\nu-\alpha \kappa)&\int_0^t \|\hat{\bzeta}\|^2ds\leq C(\kappa,\nu,\alpha)h^2\displaystyle{\int_0^t}\big(\| \nabla\hat \bzeta\|^2
+h^2\|\nabla\hat p\|^2+\kappa h^2\,\|\tilde{\Delta}\hat	\bu\|^2+\kappa h^2\,\|e^{\alpha s}\tilde{\Delta} \bu_s\|^2\big)ds.
\end{align}
Since $0\leq \alpha < \displaystyle{\frac{\nu \lambda_1}{4\big(1+\kappa \lambda_1\big)}}$, $ (\nu-\alpha \kappa) > 0.$ Then, 
use 
estimates from Lemmas \ref{l1} and \ref{L61} to complete
the rest of the proof.\hfill{$\Box$}
\begin{remark}\label{lr2}
 Here again, using (\ref{ea14})-(\ref{ea16}), we tackle the first term on the right-hand side of (\ref{ea19*}) and express 
 the entire right-hand side of (\ref{ea19*}) as an integration. As mentioned earlier in Remark \ref{lr1}, this 
 provides $\bs(t)$ in the estimate which is used to handle regularity issues of the solution at $t=0$ in the process of making 
 error bounds independent of $\kappa$.
\end{remark}
\noindent
Now, Lemma \ref{L63*} provides the estimate for the time derivative $\bzeta_t$. Here again, the estimate differs from the 
estimate of $\bzeta_t$ in Lemma 5.3 of \cite{ANS6} in terms of involvement of weight function $\bs$ and additional power 
of 
$\kappa$. As stated earlier, the 
presence of $\bs(t)$ and additional power of $\kappa$ in the estimate play a crucial role in making error estimates 
independent of $\kappa$. 
The proof proceeds in an exactly similar manner as the proof of Lemma \ref{L62*} with the right-hand side of (\ref{E624}) 
replaced by $e^{\alpha t} \bzeta_t$. But in order to justify the presence of $\bs$ in the estimate, we present a short proof.\\
 \begin{ldf}\label{L63*}
Under the assumptions {\rm{({\bf {A1}})}}-{\rm{({\bf {A2}})}}  and {\rm{({\bf B1})-({\bf B2})}}, there exists a positive constant 
$C=C(\kappa,\nu,\alpha,\lambda_1, M)$ such that for $\displaystyle{0\leq\alpha< \frac{\nu\lambda_1}{4\big(1+\kappa\lambda_1\big)}}$, 
the following holds true:
\begin{eqnarray}
\bs^{-1}(t)\kappa^2\int_0^t e^{2\alpha s}\|(\bu_s-V_h\bu_s)(s)\|^2 ds \leq C h^4.\nonumber
\end{eqnarray}
\end{ldf}
\noindent
 {\bf Proof.} For obtaining the desired estimate of $\bzeta_t$, we replace the right-hand side of (\ref{E624}) by 
 $e^{\alpha t} \bzeta_t$ and form an $\bL^2$-inner product of resulting equation with $e^{\alpha t}\bzeta_t$. Then, use 
 (\ref{E616}) with ${\bphi}_h= e^{\alpha t}P_h \bw$ in a similar way as in the $\bL^2$-estimate of $\bzeta$ to obtain
 \begin{align}\label{E644}
\|e^{\alpha t}\bzeta_t\|^2 &=\nu a(\bw-P_h \bw,e^{\alpha t}\bzeta_t)-(q-j_h q,,e^{\alpha t}\nabla\cdot \bzeta_t)-\frac{\nu^2}{\kappa}\,a(e^{\alpha t}\bzeta,P_h\bw)+\frac{\nu}{\kappa}(\hat p,\nabla\cdot P_h \bw)\nonumber\\
&=\nu a(\bw-P_h \bw,e^{\alpha t}\bzeta_t)-(q-j_h q,e^{\alpha t}\nabla \cdot\bzeta_t)-\frac{\nu^2}{\kappa}\,a(e^{\alpha t}\bzeta,P_h\bw-\bw)
-\frac{\nu^2}{\kappa}\,a(e^{\alpha t}\bzeta,\bw)\nonumber\\
&+\frac{\nu}{\kappa}(\hat p-j_h \hat p, \nabla\cdot (P_h\bw-\bw)).
\end{align}
Multiply (\ref{E644}) by $\kappa$ and use Cauchy-Schwarz's inequality with (\ref{3.5}), approximation property 
({\bf B1}), regularity result (\ref{E6281*}) with right-hand side $e^{\alpha t}\bzeta_t$. Then, after squaring 
both sides of the resulting equation, perform an integration with respect to time from $0$ to $t$ to obtain 
 \begin{align}
\label{E645e}
\kappa^2\int_0^t \|e^{\alpha s}\bzeta_s(s)\|^2ds& \leq C(\nu)\int_0^t 
 h^2\big(\kappa^2\|e^{\alpha s}\nabla\bzeta_s(s)\|^2+ \| \nabla\hat \bzeta(s)\|^2+\|\hat \bzeta(s)\|^2+h^2\|\nabla\hat p(s)\|^2 
\big)ds.
\end{align}
An application of Lemmas \ref{l1}, \ref{L61}, \ref{L62}, \ref{L62*} would lead to the desired estimates. \hfill{$\Box$}

\noindent\begin{ldf}\label{La2}
	Under the assumptions {\rm{({\bf {A1}})}}-{\rm{({\bf {A2}})}}  and {\rm{({\bf B1})-({\bf B2})}}, there exists a positive constant 
	$C=C(\kappa,\nu,\alpha,\lambda_1, M)$ such that for $\displaystyle{0\leq\alpha< \frac{\nu\lambda_1}{4\big(1+\kappa\lambda_1\big)}}$, 
	the following estimate holds true for $t >0$:
	\begin{align}
\kappa^2\|(\bu_t-V_h\bu_t)(t)\|^2+\bs^{-1}(t)\kappa\int_0^t e^{2\alpha s}\|(\bu_s-V_h\bu_s)(s)\|^2 ds
 &\leq C h^4,\nonumber\\
 \|(\bu-V_h\bu)(t)\|^2&\leq C h^4.\nonumber
	\end{align}
\end{ldf}
\noindent{\it Proof.} 
Once again, we apply the Aubin-Nitsche duality
argument. Let $(\bw,q)$ be the unique
solution of the following steady state Stokes system:
\begin{align}
-\nu\Delta \bw+\nabla q&=\bzeta_t\ \ {\rm in}\ \ \Omega, \label{ea25}\\
\nabla \cdot \bw&=0\ \ \ \ {\rm in}\ \ \Omega,  \label{ea26} \\
\bw|_{\partial\Omega}&=0.  \label{ea27}
\end{align}
Now, using assumption ({\bf A1}), $(\bw,q)$ satisfies the following
regularity result:
\begin{eqnarray}\label{ea28}
\|\bw\|_2+\|q\|_{H^1/{\rm I\!R}}\le C\|\bzeta_t\|.
\end{eqnarray}
Taking an $L^2$-inner product between (\ref{ea25}) and $\bzeta_t$ and using the discrete incompressibility condition, we obtain
\begin{align}\label{ea29}
\|\bzeta_t\|^2= \nu\, a (\bw-P_h \bw, \bzeta_t)-(q-j_{h}q,\nabla\cdot\bzeta_t)+\nu\, a( P_h \bw, \bzeta_t).
\end{align}
Now, by using (\ref{E616}) with ${\bphi}_h$ replaced by $P_h\bw$ and (\ref{ea26}), the last term in (\ref{ea29}) can be rewritten as 
\begin{eqnarray}\label{ea30}
\nu\, a(  \bzeta_t,P_h \bw)=(p_t-j_h p_t,\nabla \cdot (P_h \bw-\bw))-\kappa\,a( \bzeta_{tt}, P_h \bw-\bw)-\kappa\,a( \bzeta_{tt},\bw).
\end{eqnarray}
Use (\ref{ea25}) to rewrite last term in (\ref{ea30}) as 
\begin{align}\label{ee3}
 \kappa\,a(\bzeta_{tt},\bw)=\frac{\kappa}{\nu}(\bzeta_t,\bzeta_{tt})+\frac{\kappa}{\nu}(q-j_h q,\nabla\cdot\bzeta_{tt})
 =\frac{\kappa}{2\nu}\frac{d}{dt}\|\bzeta_t\|^2+\frac{\kappa}{\nu}(q-j_h q,\nabla\cdot\bzeta_{tt}).
\end{align}
Apply (\ref{ea30}), (\ref{ee3}) in (\ref{ea29}) to obtain
\begin{align}\label{ea31}
\|\bzeta_t\|^2&=\nu\, a (\bw-P_h \bw, \bzeta_t)-(q-j_{h}q,\nabla\cdot\bzeta_t)+(p_t-j_h p_t,\nabla \cdot (P_h \bw-\bw))-
\kappa\,a( \bzeta_{tt}, P_h \bw-\bw)\nonumber\\
&-\frac{\kappa}{2\nu}\frac{d}{dt}\|\bzeta_t\|^2-\frac{\kappa}{\nu}(q-j_h q,\nabla\cdot\bzeta_{tt}).
\end{align}
A simplification of (\ref{ea31}) yields
\begin{align}\label{ea32}
\|\bzeta_t\|^2+\frac{\kappa}{2\nu}\frac{d}{dt}\|\bzeta_t\|^2&=\nu\, a (\bw-P_h \bw, \bzeta_t)-(q-j_{h}q,\nabla\cdot\bzeta_t)+(p_t-j_h p_t,\nabla \cdot (P_h \bw-\bw))\nonumber\\
&-\kappa\,a( \bzeta_{tt}, P_h \bw-\bw)-\frac{\kappa}{\nu}(q-j_h q,\nabla\cdot\bzeta_{tt}).
\end{align}
After multiplying (\ref{ea32}) by $\kappa\,\bs$, rewrite the resulting equation as
\begin{align}\label{ea33}
\bs&\|\bzeta_t\|^2+\frac{\kappa}{2\nu}\frac{d}{dt}\bs\|\bzeta_t\|^2
=\frac{\kappa}{2\nu}\bs_t\|\bzeta_t\|^2+\nu\,\frac{d}{dt} \bs\,a (\bw-P_h \bw, \bzeta)-\nu\,
 \bs_t\,a (\bw-P_h \bw, \bzeta)\nonumber\\
 &-\frac{d}{dt}\bs\,(q-j_{h}q,\nabla\cdot\bzeta)
+\bs_t\,(q-j_{h}q,\nabla\cdot\bzeta)+\frac{d}{dt}\bs\,(p-j_h p,\nabla \cdot (P_h \bw-\bw))\nonumber\\
&-\bs_t\,(p-j_h p,\nabla \cdot (P_h \bw-\bw))-\kappa\,\frac{d}{dt}\bs a( \bzeta_t, P_h \bw-\bw)+\kappa\,\bs_t a( \bzeta_{t}, P_h \bw-\bw)\nonumber\\
&-\frac{\kappa}{\nu}\frac{d}{dt}\bs(q-j_h q,\nabla\cdot\bzeta_t)+\frac{\kappa}{\nu}\bs_t(q-j_h q,\nabla\cdot\bzeta_t).
\end{align}
An integration of (\ref{ea33}) with respect to time from $0$ to $t$ along with a use of Cauchy-Schwarz's inequality, 
Young's inequality leads to
\begin{align}\label{ea35}
\kappa&\displaystyle{\int_0^t}\bs\|\bzeta_s\|^2ds +\frac{\kappa^2}{2\nu}\bs\|\bzeta_t\|^2
\leq \frac{\kappa^2}{2\nu}\displaystyle{\int_0^t}\left(\bs_s\|\bzeta_s\|^2\right)ds+\nu\,\kappa\,\bs h\|\tilde\Delta\bw\|\|\nabla\bzeta\|
+\kappa \bs h\|\nabla q\|\|\nabla\cdot\bzeta\|\nonumber\\
&+\kappa \bs h^2\|\nabla p\|\|\tilde{{\Delta}}\bw\|+\kappa^2\,h^2\,\bs \|\nabla \bzeta_t\|\|\tilde\Delta\bw\|
+C(\nu)\,\kappa^2\bs\|\nabla q\|\|\nabla\cdot\bzeta_t\|
+h^2\displaystyle{\int_0^t}\bs_s\|\nabla \bzeta\|^2ds\nonumber\\
&+\kappa^2\displaystyle{\int_0^t}\bs_s\|\nabla q\|^2ds+h^4\displaystyle{\int_0^t}\bs_s\|\nabla p\|^2ds
+\kappa^2\displaystyle{\int_0^t}\bs_s\|\tilde{\Delta}\bw\|^2ds+\kappa^2\displaystyle{\int_0^t}\bs_s
\|\nabla\bzeta_s\|^2ds.
\end{align}
Apply the bounds from Lemmas \ref{l1}, \ref{L61}, \ref{L62}, \ref{L62*}, \ref{L63*} and 
the regularity estimates (\ref{ea28}) to arrive at the desired result.\\
\noindent
Next, to derive $\L^{\infty}(\bL^2)$ estimates of $\bzeta$, follow the similar steps as in Lemma \ref{L62*} 
with $\hat\bzeta$ replaced by $\bzeta$ and arrive 
 at (\ref{E62810*}). Then, use  Cauchy-Schwarz's inequality to obtain
 \begin{align}\label{ea38}
  \|\bzeta\|^2&\leq  Ch\left(\nu\| \nabla\bzeta\|+h\|\nabla p\|+\kappa\,\|\nabla\bzeta_t\|+
  \kappa\|\bzeta_t\|\right)(\|\bw\|_2+\|q\|_1).
  \end{align}
  Apply regularity estimates (\ref{E6281*}) with right-hand side as $\bzeta$ along with (\ref{ea13}), Lemmas \ref{l1}, 
  \ref{L62}, 
  \ref{La2} to compelte the rest part of the proof.\hfill{$\Box$}

  \begin{ldf}\label{La4}
  	Under the assumptions {\rm{({\bf {A1}})}}-{\rm{({\bf {A2}})}}  and {\rm{({\bf B1})-({\bf B2})}}, there exists a positive constant 
  	$C=C(\kappa,\nu,\alpha,\lambda_1,M)$ such that for $\displaystyle{0\leq\alpha< \frac{\nu\lambda_1}{4\big(1+\kappa\lambda_1\big)}}$, 
  	the following holds true:
  	\begin{eqnarray}
  	\bs^{-1}(t)\int_0^t \bs_1(s)\|(\bu_s-V_h\bu_s)(s)\|^2 ds \leq C h^4,\nonumber
  	\end{eqnarray}
  	where $\bs_1(t):=\tau^2(t)e^{2\alpha t}$ with $\tau(t):=\min\{t,\,1\}$. 
  \end{ldf}
\noindent{\it Proof.} 
A use of Cauchy-Schwaz's inequality,  {\rm{({\bf B1}) and (\ref{3.5}) in (\ref{ea31}) yield
\begin{align}\label{ea46}
\|\bzeta_t\|^2&\leq \nu\,h\|\tilde{\Delta}\bw\|\|\nabla\bzeta_t\|+h\|\nabla q\|\|\nabla\cdot\bzeta_t\|+h^2
\|\nabla p_t\|\|\tilde{\Delta}\bw\|\nonumber\\
&+h\kappa\,\|\nabla\bzeta_{tt}\|\|\tilde{\Delta}\bw\|
-\frac{\kappa}{2\nu}\frac{d}{dt}\|\bzeta_t\|^2+\frac{\kappa}{\nu}h\|\nabla q\|\|\nabla\cdot\bzeta_{tt}\|.
\end{align}
Multiply (\ref{ea46}) by $\bs_1$, apply (\ref{ea42}) and integrate the resulting equation from $0$ to $t$  to arrive at
\begin{align}\label{ea47}
\displaystyle{\int_0^t}&\bs_1\|\bzeta_s\|^2ds\leq \nu\,h\displaystyle{\int_0^t}\bs_1\|\tilde{\Delta}\bw\|\|\nabla\bzeta_s\|ds
+h\displaystyle{\int_0^t}\bs_1\|\nabla q\|\|\nabla\cdot\bzeta_s\|ds\nonumber\\
&+h^2\displaystyle{\int_0^t}\bs_1\|\nabla p_s\|\|\tilde{\Delta}\bw\|ds
+h\,\displaystyle{\int_0^t}\bs_1(\nu\|\nabla\bzeta_{s}\|+\|\nabla p_s\|+\kappa \|\tilde{\Delta }\bu_{ss}\|)\|\tilde{\Delta}\bw\|ds\nonumber\\
&-\frac{\kappa}{2\nu}\bs_1\|\bzeta_s\|^2+\frac{\kappa}{2\nu}\displaystyle{\int_0^t}\bs_{1,s}\|\bzeta_s\|^2ds
+C(\nu)h\displaystyle{\int_0^t}\bs_1 (\nu\|\nabla\bzeta_{s}\|+\|\nabla p_s\|+\kappa \|\tilde{\Delta }\bu_{ss}\|)
\|\nabla q\|ds.
\end{align}
An application of Young's inequality along with regularity estimates (\ref{ea28}) leads to
\begin{align}\label{ea48}
\displaystyle{\int_0^t}&\bs_1\|\bzeta_s\|^2ds\leq  C(\nu)\left( \displaystyle{\int_0^t}
\bs_1(h^2\|\nabla\bzeta_s\|^2+h^4\|\nabla p_s\|+ h^4\kappa\|\tilde{\Delta }\bu_{ss}\|^2)+\kappa \bs_{1,s}\|\bzeta_s\|^2\right)ds.
\end{align}
A use of Lemmas \ref{l2}, \ref{La3}, \ref{La2} would lead us to the desired result. \hfill{$\Box$}\\

\noindent 
 Since $\bxi=\bzeta + \brho$ and the estimates of $\bzeta$ are already 
derived, it suffices to derive the estimates of $\brho$ to obtain estimates for $\bxi$. Below, in Lemma \ref{L65*}, we state 
without proof estimates of $\brho$.  We skip the proof as it follows the similar lines as in the proofs of 
Lemma 5.6 (\cite{bnpdy}) and Lemma \ref{L61} in this article. We also present a couple of estimates of $\bxi$ which can 
be easily derived using the estimates of $\bzeta$ and $\brho$.   
\begin{ldf}\label{L65*}
Under the assumptions {\rm{({\bf {A1}})}}-{\rm{({\bf {A2}})}}  and {\rm{({\bf B1})-({\bf B2})}}, there exists a positive constant 
$C=C(\kappa,\nu,\alpha,\lambda_1,M)$ such that for $\displaystyle{0\leq\alpha< \frac{\nu}{4\big(1+\kappa\lambda_1\big)}}$, 
the following estimates hold true:
\begin{align*}
 \kappa^2(\|\brho(t)\|^2+\kappa\|\nabla \brho(t)\|^2)+\kappa^2\,e^{-2\alpha t}\displaystyle{\int_0^t}e^{2\alpha s}\|\nabla  \brho(s)\|^2 ds&
 \leq C h^4,\\
\bs^{-1}(t)\kappa^2 \displaystyle{\int_0^t}e^{2\alpha s}\|\bxi_s(s)\|^2ds &\leq C\,h^4,\\ 
\bs^{-1}(t)\displaystyle{\int_0^t}e^{2\alpha s}\left(\kappa\,\|\bxi_s(s)\|^2+\kappa^2\,\|\nabla\bxi_s(s)\|^2+\|\nabla\bxi(s)\|^2\right)ds&\leq 
C\,h^2.
\end{align*}
\end{ldf}
\hfill{$\Box$}
\begin{ldf}\label{L65n}
Under the assumptions {\rm{({\bf {A1}})}}-{\rm{({\bf {A2}})}}  and {\rm{({\bf B1})-({\bf B2})}}, there exists a positive constant 
$C=C(\kappa,\nu,\alpha,\lambda_1,M)$ such that for $\displaystyle{0\leq\alpha< \frac{\nu}{4\big(1+\kappa\lambda_1\big)}}$, 
the following estimate holds true:
\begin{align*}
\bs^{-1}(t)\displaystyle{\int_0^t}e^{2\alpha s}\,\|\bxi(s)\|^2ds&\leq C\,h^4.
 \end{align*}
\end{ldf}
\noindent
\noindent{\it Proof.}
To estimate $L^2$-error, we use the following duality argument: For fixed $t>0$ with $t\in(0,T)$, 
let $\bw(\tau)\in\bJ_1$, $q(\tau)\in L^2 / {\rm I\!R}$ be the unique solution of the 
backward Stokes problem
\begin{align}\label{ea56}
\bw_{\tau}+\nu\,\tilde{\Delta}\bw-\nabla q=e^{2\alpha t}\bxi,\,0\leq \tau \leq t,\,\,\bw(t)=0.
\end{align}
The pair $(\bw,q)$ satisfy the following regularity estimates
\begin{align}\label{ea57}
\displaystyle{\int_0^t}e^{-2\alpha \tau}(\|\tilde{\Delta}\bw\|^2+\|\bw_\tau\|^2+\|\nabla q\|^2)d\tau\leq C
\displaystyle{\int_0^t}e^{2\alpha \tau}\|\bxi\|^2d\tau.
\end{align}
Form an $L^2$ inner product between (\ref{ea56}) and $\bxi$ to arrive at
\begin{align}\label{ea58}
e^{2\alpha \tau}\|\bxi\|^2&=(\bxi,\bw_{\tau})-\nu\,a(\bxi,\bw)+(q,\nabla\cdot \bxi)\nonumber\\
&=\frac{d}{d\tau}(\bxi,\bw)-(\bxi_\tau,\bw-P_h\bw)-\nu\,a(\bxi,\bw-P_h\bw)\nonumber\\
&+(q-j_hq,\nabla\cdot \bxi)-(\bxi_\tau,P_h\bw)-\nu\,a(\bxi,P_h\bw).
\end{align}
A use of (\ref{E615}) with $\bphi_h$ replaced by $P_h\bw$ in (\ref{ea58}) yields
\begin{align}\label{ea59}
e^{2\alpha \tau}\|\bxi\|^2&=\frac{d}{d\tau}(\bxi,\bw)-(\bxi_{\tau},\bw-P_h\bw)-\nu\,a(\bxi,\bw-P_h\bw)+(q-j_hq,\nabla\cdot \bxi)\nonumber\\
&+\kappa\,a(\bxi_{\tau},P_h\bw-\bw)-(p-j_hp,\nabla\cdot (P_h\bw-\bw))+\kappa\,a(\bxi_{\tau},\bw).
\end{align}
Note that, 
\begin{align}\label{ea60}
(\bxi_{\tau},\bw-P_h\bw)=\frac{d}{d\tau}(\bxi,\bw-P_h\bw)-(\bu-P_h\bu,\bw_{\tau}).
\end{align}
A simplification of (\ref{ea59}), using (\ref{ea60}) leads to
\begin{align}\label{ea61}
e^{2\alpha \tau}\|\bxi\|^2&=\frac{d}{d\tau}(\bxi,P_h\bw)+(\bu-P_h\bu,\bw_{\tau})-\nu\,a(\bxi,\bw-P_h\bw)+(q-j_hq,\nabla\cdot \bxi)\nonumber\\
&+\kappa\,a(\bxi_\tau,P_h\bw-\bw)-(p-j_hp,\nabla\cdot (P_h\bw-\bw))-\kappa\,(\bxi_\tau,\tilde{\Delta}\bw).
\end{align}
An integration of (\ref{ea61}) with respect to time from $0$ to $t$ along with Cauchy-Schwarz's inequality yields
\begin{align}\label{ea62}
\displaystyle{\int_0^t}e^{2\alpha \tau}\|\bxi(\tau)\|^2ds\leq&(\bxi(t),P_h\bw(t))-(\bxi(0),P_h\bw(0))+\displaystyle{\int_0^t}(h^2\|\tilde{\Delta}\bu\|
 +\nu\,h\|\nabla \bxi\|+h\|\nabla\cdot \bxi\|\nonumber\\
 &+\kappa\,h\|\nabla\bxi_\tau\|+h^2\|\nabla\,p\|+\kappa\,\|\bxi_\tau\|)(\|\tilde{\Delta}\bw\|+\|\nabla\,q\|)d\tau. 
\end{align}
The first term in (\ref{ea62}) vanishes due to $\bw(t)=0$ and the second term disappears due to the orthogonality property of $P_h$. 
Now, a use of Young's inequality along with the regularity estimates (\ref{ea57}) leads to
\begin{align}\label{ea63}
\displaystyle{\int_0^t}e^{2\alpha \tau}\|\bxi(\tau)\|^2d\tau\leq&C(\nu)\bigg(\displaystyle{\int_0^t}(h^4\|\tilde{\Delta}\hat\bu\|^2
+\nu\,h^2\,\|\nabla \hat\bxi\|^2+h^2\|\nabla\cdot\hat \bxi\|^2\nonumber\\
&+\kappa^2\,h^2\,e^{2\alpha \tau}\|\nabla\bxi_\tau\|^2+h^4\|\nabla\,\hat p\|^2+\kappa^2\,e^{2\alpha s}\|\bxi_\tau\|^2\bigg)d\tau. 
\end{align}
Apply estimates from Lemmas \ref{l1}, \ref{L65*} to arrive at the desired result.\hfill{$\Box$}\\
\begin{ldf}\label{L63}
	Let the assumptions {\rm{({\bf {A1}})}}-{\rm{({\bf {A2}})}}  and {\rm{({\bf B1})-({\bf B2})}} be satisfied. Then, 
	there exists a positive constant $C=C(\kappa,\nu,\alpha,\lambda_1,M)$ such that 
	for $\displaystyle{0\leq\alpha< \frac{\nu\lambda_1}{4\big(1+\kappa \lambda_1\big)}}$, the following estimates hold true:
	\begin{eqnarray}
	&&\kappa(\|\brho(t)\|^2+\kappa\|\nabla \brho(t)\|^2)+\kappa\,e^{-2\alpha t}\displaystyle{\int_0^t}
	e^{2\alpha s}\|\nabla\brho(s)\|^2 ds
\leq Ch^4,\label{1n}\\
	&&\|\brho(t)\|^2+\kappa\,\|\nabla \brho(t)\|^2+ \nu\,\bs_1^{-1}(t) \displaystyle{\int_0^t }\bs_1(s)\|\nabla\brho(s)\|^2ds
	\leq C h^4,\label{2n}
	\end{eqnarray}
	where $\bs_1(t):=\tau^2(t)e^{2\alpha t}$ with $\tau(t):=\min\{t,\,1\}$. 
\end{ldf}
\noindent{\it Proof.} 
Subtracting (\ref{E616}) from (\ref{E615}), we find that
\begin{eqnarray}
\label{E647}
(\brho_t,\bphi_h)+\kappa\,a(\brho_t, {\bphi}_h)+ \nu a ( \brho, {\bphi}_h) = -(\bzeta_t, {\bphi}_h) \;\;\;\;\; \forall {\bphi}_h\in \bJ_h.
\end{eqnarray}
Multiply (\ref{E647}) by $\kappa\,e^{2\alpha t}$, substitute $\bphi_h=\brho$ and use Cauchy-Schwarz's inequality, 
Young's inequality in the resulting equation. 
Then, integrate 
the equation from $0$ to $t$ to arrive at
\begin{align}
\label{ea52}
\kappa(\|\hat \brho\|^2+\kappa\|\nabla\hat \brho\|^2)+\kappa\displaystyle{\int_0^t}\|\nabla \hat \brho\|^2 ds&\leq 
C(\kappa,\alpha,\lambda_1)\,\int_0^t\,\left(\kappa^2\| e^{\alpha s}\bzeta_s(s)\|^2+
\| \hat \brho\|^2 \right)ds.
\end{align}
Note that, $\brho=\bxi-\bzeta$. A use of the triangle inequality along with Lemmas \ref{L62*} and \ref{L65n} yields
\begin{align}\label{ee5}
 \displaystyle{\int_0^t}\,\| \hat \brho\|^2 ds\leq \displaystyle{\int_0^t}\,(\| \hat \bzeta\|^2+\|\hat\bxi\|^2) 
 ds\leq C h^4\,\bs.  
\end{align}
An application of the results from Lemma \ref{L63*} and (\ref{ee5}) in (\ref{ea52}) and a 
multiplication of the resulting equation by $e^{-2\alpha t}$ complete the proof of (\ref{1n}).\\
\noindent
Next to prove (\ref{2n}), substitute $\bphi_h=\brho$ in (\ref{E647}) and multiply the resulting equation by 
$\bs_1$ to arrive at
\begin{align}
\frac{1}{2}\frac{d}{dt}\bs_1(\|\brho\|^2+\kappa\,\|\nabla \brho\|^2)+ \nu \bs_1\|\nabla\brho\|^2
= -\bs_{1}(\bzeta_t, \brho) +\bs_{1,t}(\|\brho\|^2+\kappa\,\|\nabla \brho\|^2).
\end{align}
After applying Cauchy-Schwarz's inequality and Young's inequality, integrate the resulting equation with respect 
to time from $0$ to $t$ to obtain
\begin{align}
\bs_1(\|\brho\|^2+\kappa\,\|\nabla \brho\|^2)+ \nu \displaystyle{\int_0^t }\bs_1\|\nabla\brho\|^2ds
\leq C\displaystyle{\int_0^t}\left(\frac{\bs^2_{1}}{\bs_{1,s}}\|\bzeta_s\|^2 +\bs_{1,s}(\|\brho\|^2
+\kappa\,\|\nabla \brho\|^2)\right)ds.\nonumber
\end{align}
Use estimates from (\ref{1n}), (\ref{ee5}) and Lemma \ref{La4} to arrive at (\ref{2n}) and this completes the proof of 
Lemma \ref{L63}.\hfill{$\Box$}\\
\noindent
Now, we derive the proof of the main Theorem \ref{T31*}.

\noindent
Note that $\e=\bu-\bu_h=(\bu-\bvh)+(\bvh-\bu_h)=\bxi+\bta$. A use of the triangle inequality, the inverse inequality and 
Lemmas \ref{La2}, \ref{L63} 
lead to the following estimates of $\bxi$.
\begin{eqnarray}\label{an1}
 \|\bxi(t)\|^2+h^2\,\|\nabla\bxi(t)\|^2\leq Ch^4.
\end{eqnarray}
In Lemma \ref{L63n}, we present the estimates of $\eta$. For a proof, one may refer to 
\cite{ANS6} 
(Theorem 5.1, pp. 249 - 250).
\begin{ldf}\label{L63n}
	Let the assumptions {\rm{({\bf {A1}})}}-{\rm{({\bf {A2}})}}  and {\rm{({\bf B1})-({\bf B2})}} be satisfied. Then, 
	there exists a positive constant $C=C(\kappa,\nu,\alpha,\lambda_1,M)$ such that 
	for $\displaystyle{0\leq\alpha< \frac{\nu\lambda_1}{4\big(1+\kappa \lambda_1\big)}}$, the following estimates hold true:
	\begin{eqnarray}
	&&\|\eta(t)\|^2+\kappa\|\nabla \eta(t)\|^2+\kappa\,e^{-2\alpha t}\displaystyle{\int_0^t}
	e^{2\alpha s}\|\nabla\eta(s)\|^2 ds
\leq K(t)h^4.\nonumber
	\end{eqnarray}
\end{ldf}
\noindent
Moreover, under the assumptions of Theorem \ref{T31*} and the uniqueness condition (\ref{tn}), 
the constant $K(t)=C$. That is, the estimates are valid uniformly with respect to time. \\

\noindent{\it Proof of Theorem \ref{T31*}}. The proof follows by using the triangle inequality, inverse inequality, 
(\ref{an1}) and Lemma \ref{L63n}.\hfill{$\Box$}\\

\noindent
Following the similar steps as in \cite{ANS6} (Theorem 6.1) and using $\kappa$ independent estimates derived earlier, 
we arrive at the desired pressure error estimates in Theorem \ref{T32*} and this completes the proof. \hfill{$\Box$}  
  \section{ Fully Discrete Approximation}
\setcounter{tdf}{0}
\setcounter{ldf}{0}
\setcounter{cdf}{0}
\setcounter{equation}{0} 
In this section, we apply a backward Euler method for time discretization of the finite element Galerkin 
approximation (\ref{3.1}) of (\ref{1.1})-(\ref{1.3}). Let $\{t_n\}_{n=0}^{N}$ be a uniform partition of $[0,T]$, and $t_n=nk$, with 
time step $k>0$. For smooth function $\bphi$ defined on $[0, T]$, set $\bphi^n=\bphi(t_n)$ and 
$\bar \partial _t \bphi^n=\frac{(\bphi^n-\bphi^{n-1})}{k}$.\\
\noindent
The backward Euler method applied 
to (\ref{3.1}) determines a sequence of functions  ${\{\bU^n\}}_{n\geq1}\in \bH_h$ 
and ${\{P^n\}}_{n\geq1}\in L_h$ as solutions of the following recursive nonlinear algebraic equations: 
 \begin{align}\label{4.2a1}
(\bar\partial_t {\bU}^n,\bphi_h)&+\kappa a (\bar\partial_t {\bU}^n,\bphi_h)+\nu a({\bU}^n,\bphi_h)\nonumber\\
&+b({\bU}^n, {\bU}^n,\bphi_h)=(P^n,\nabla\cdot \bphi_h)+(\f^n,\bphi_h) \,\,\,\,\forall \bphi_h\in \bH_h,\\
& (\nabla \cdot\bU^n, \chi_h)=0\,\,\,\,\,\,\,\, \,\,\,\,\,\,\,\forall \chi_h\in L_h,\nonumber \\
 &\bU^{0}=\bu_{0h}.\nonumber
\end{align}
Equivalently, we seek ${\{\bU^n\}}_{n\geq1}\in \bJ_h$ such that
\begin{align}\label{4.2a}
(\bar\partial_t {\bU}^n,\bphi_h)+\kappa a (\bar\partial_t {\bU}^n,\bphi_h)+\nu a({\bU}^n,\bphi_h)&
+b({\bU}^n, {\bU}^n,\bphi_h)=(\f^n,\bphi_h) \,\,\,\,\forall \bphi_h\in \bJ_h, \\
\bU^{0}&=\bu_{0h}.\nonumber
\end{align}
Next, in Lemma \ref{ll1}, we state a {\it priori} bounds for the discrete solution $\{\bU^n\}_{n\geq 1}$. We skip the proof 
as it will be an imitation of the proof of Lemma 4.1 in \cite{P4}.
\begin{ldf}\label{ll1}
With $\displaystyle{0\le\alpha< \frac{\nu \lambda_1}{4(1+\lambda_1\kappa)}}$, choose $k_0$ so that 
for $0<k\leq k_0$
\begin{align}
\label{5.1}
\displaystyle{\frac{\nu k\lambda_1}{\kappa \lambda_1+1}+1> e^{\alpha k}}.
 \end{align}
 Then the discrete solution $\bU^N$, $N \geq 1$ of (\ref{4.2a}) satisfies
\begin{eqnarray}\label{1a}
(\|{\bU^N}\|^2+\kappa\|\nabla{\bU^N}\|^2)+ e^{-2\alpha t_{N}}~k \displaystyle{\sum_{n=1}^{N}} 
e^{2\alpha t_n}\|\nabla {\bU}^n\|^2\leq C(\nu,\alpha,\lambda_1)e^{-2\alpha t_N}(\|{\bU}^0\|^2+\kappa 
\|\nabla {\bU}^0\|^2+\|\f\|^2_{\infty}).\nonumber
\end{eqnarray}
\begin{align*} 
\end{align*}
\hfill{$\Box$}
\end{ldf}
\noindent
Next, we proceed to derive fully discrete estimates for the velocity error $\e^n=\bU^n-\bu_h(t_n)=\bU^n-\bu_h^n$ 
and for the pressure error $\rho^n=P^n-p_h(t_n)=P^n-p_h^n$. Below, in Lemma \ref{T51}, we present the 
various estimates of $\e^n$. The proof of (\ref{en5}) follows the similar lines as in the proof of Theorem 5.1 of \cite{P4}. Therefore, we skip the proof. 
 The estimates of $\|\bar\partial_t\e^n\|$ and $\kappa\|\bar\partial_t\nabla\e^n\|$ are also discussed in Lemma 5.1 of
\cite{P4}, but, these estimates involve $\kappa^{-1}$ term. Therefore, here we provide a short proof of (\ref{en6}) 
by only highlighting the steps involved in making estimates independent of the inverse power of $\kappa$. 
\begin{ldf}\label{T51}
Let $\displaystyle{0\leq\alpha< \frac{\nu \lambda_1}{4(1+\kappa\lambda_1)}}$ and $k_0>0$ be such that for $0<k\leq k_0$, 
(\ref{5.1}) is satisfied.
For some fixed $h>0$, let $\bu_h(t)$ satisfies (\ref{3.2}).  Then, there is a positive constant $C_T$ that depends on $T$ 
such that
 \begin{eqnarray}
  \|\e^i\|^2+\kappa \|\nabla\e^i\|^2+ke^{-2\alpha t_n}\displaystyle{\sum_{i=1}^n}e^{2\alpha t_i}\|\nabla\e^i\|^2&\leq & 
  C_T k^2,\label{en5}\\
  \|\bar\partial_t\e^n\|^2_{-1}+\kappa^2\,\|\bar\partial_t\nabla {\e}^n\|^2&\leq & C_Tk.\label{en6}
 \end{eqnarray}
\end{ldf}
\noindent{\it Proof.} To prove (\ref{en6}), consider (\ref{3.2}) at $t=t_i$ and subtract it from (\ref{4.2a}) to obtain 
\begin{align}\label{en1}
 (\bar \partial_t\e^n,\bphi_h)&+\kappa a(\bar \partial_t\e^n,\bphi_h)+\nu a(\e^n,\bphi_h)\\
&=(\bs_1^n,\bphi_h)+\kappa a(\bs_1^n,\bphi_h)+\Lambda_h(\bphi_h)\,\, \forall \bphi_h \in \bJ_h,\nonumber
\end{align}
where $\bs_1^n= \bu_{ht}^{n}-\bar\partial_t \bu_h^n$ and $\Lambda_h(\bphi_h)=
b(\bu_h^n,\bu_h^n,\bphi_h)-b(\bU^n,\bU^n,\bphi_h)$.\\
Note that, applying Taylor's series expansion in the interval $(t_{i-1},t_i)$, Cauchy-Schwarz's inequality, 
Young's inequality and estimates from Lemma \ref{L72}, we arrive at
\begin{eqnarray}\label{en2}
|(\bs_1^{n},\bphi_h)|&\leq& \frac{1}{k}\displaystyle{\int_{t_{n-1}}^{t_n}}(t-t_{n-1})\|\bu_{htt}\|_{-1}dt\|\nabla\bphi_h\|\nonumber\\
 & \leq & C\,k^{1/2}\left\{ \displaystyle{\int_{t_{n-1}}^{t_n}}\|\bu_{htt}\|^2_{-1}dt\right\}^{1/2}\|\nabla\bphi_h\|
 \leq C\,k^{1/2}\|\nabla\bphi_h\| 
\end{eqnarray}
and 
\begin{eqnarray}\label{en3}
|\kappa\,a(\bs_1^{n},\bphi_h)|&\leq& \frac{1}{k}\displaystyle{\int_{t_{n-1}}^{t_n}}(t-t_{n-1})(\kappa\|\nabla\bu_{htt}\|)dt\|
\nabla\bphi_h\|\nonumber\\
 & \leq & C\,k^{1/2}\left\{ \kappa^2\displaystyle{\int_{t_{n-1}}^{t_n}}\|\nabla\bu_{htt}\|^2dt\right\}^{1/2}\|\nabla\bphi_h\|
 \leq C\,k^{1/2}\|\nabla\bphi_h\| . 
\end{eqnarray}
Rewrite the nonlinear term and apply generalized H\"{o}lder's inequality to observe that 
\begin{align}\label{5.17}
|\Lambda_h(\bphi_h)|&= |b(\bu_h^n,\bu_h^n,\bphi_h)-b(\bU^n,\bU^n, \bphi_h)|\nonumber\\
&=|-b(\bu_h^n,\e^n,\bphi_h)-b(\e^n,\bU^n,\bphi_h)|\leq C \left(\|\nabla\bu_H^n\|+\|\nabla\bU^n\|\right)
\|\nabla\e^n\|\|\nabla\bphi_h\|.
\end{align}
Now, substitute $\bphi_h=\bar\partial_t \e^n$ in (\ref{en1}), drop the first term from left hand side and use 
(\ref{en2})-(\ref{5.17}) to observe that
 \begin{align}
 \kappa\,\|\bar\partial_t\nabla\e^n\|^2\leq  C\left(\|\nabla\e^n\|+k^{1/2}
 +(\|\nabla\bu_H^n\|+\|\nabla\bU^n\|)\|\nabla\e^n\|\right)\|\bar\partial_t\nabla\e^n\|.\nonumber 
 \end{align}
 A use of (\ref{en5}), Lemmas \ref{L72}, \ref{ll1} yield
 \begin{align}\label{e7}
 \kappa\,\|\bar\partial_t\nabla\e^n\|\leq C_Tk^{1/2}. 
 \end{align}
 Now, following the steps involved in arriving at the equation (107) from (106) in the proof Lemma 5.1 of \cite{P4}, 
 we arrive at
 \begin{align}\label{e8}
 \|\bar\partial_t\e^n\|\leq C_Tk^{1/2}.
 \end{align}
A combination of (\ref{e7}) and (\ref{e8}) completes the rest of the proof. \hfill{$\Box$}\\
\begin{remark}\label{rr3}
 Note that in the proof of Theorem 5.1 of \cite{P4}, the presence of $\kappa^{-1}$ in the first term of right 
 hand side of equation (93) is a typo, as the first term is a combination of equation (88) and (89), in which 
 estimates are independent of $\kappa^{-1}$.   
\end{remark}

\noindent
To prove the pressure error estimates, subtract (\ref{4.2a1}) from (\ref{3.1}) and write $\brho^n=P^n- p_h^n$ to obtain 
\begin{eqnarray}
 (\brho^n,\nabla\cdot\bphi_h)&= (\bar\partial_t \e^n,\bphi_h)+\kappa\,a(\bar\partial_t \e^n,\bphi_h)+\nu a(\e^n,\bphi_h)
 -\Lambda_h(\bphi_h)-(\bs_1^n,\bphi_h)-\kappa\,a(\bs_1^n,\bphi_h).\nonumber
\end{eqnarray}
A use of Cauchy-Schwarz's inequality along with (\ref{en2})-(\ref{5.17}), Lemmas \ref{L72}, \ref{ll1}, \ref{T51}  yields
\begin{eqnarray}\label{a60}
 \|\brho^n\|\leq C(\kappa,\nu,\lambda_1,M) k^{1/2}.
\end{eqnarray}
A combination of (\ref{a60}), Lemma \ref{T51} and Theorems \ref{T31*}, \ref{T32*} lead to the following 
fully discrete error estimates.
\begin{tdf}\label{Tbd1}
Under the assumptions of Theorem \ref{T31*} and Lemma \ref{T51}, the following hold true:  
\begin{eqnarray}
   &&\|\bu(t_n)-\bU^n\|\leq C (h^2+k),\qquad \|\nabla(\bu(t_n)-\bU^n)\|\leq C (h+k).\nonumber\\
   &&\|p(t_n)-P^n\|\leq C (h+k^{1/2}).\nonumber
\end{eqnarray}
\end{tdf}
\section{Numerical Experiments} 
\setcounter{tdf}{0}
\setcounter{ldf}{0}
\setcounter{cdf}{0}
\setcounter{equation}{0}
This section conducts numerical experiments to validate our theoretical results obtained in Theorem \ref{Tbd1} for finite 
element Galerkin approximations of (\ref{1.1})-(\ref{1.3}). We apply mixed finite element $P_2$-$P_0$ for space discretization 
and backward Euler method for time discretization. 
\begin{example}\label{7e1}
In this example, we choose right-hand side function $\f$ in such a way that the exact solution 
$(\bu, p) = ((u_1, u_2), p)$ takes the following form: 
\begin{align}
&u_1 = 10x^2 (x-1)^2 y(y-1)(2y -1)\,\cos t,& p = 40xy\, \cos t,\nonumber\\
&u_2 = 10y^2(y - 1)^2 x(x - 1)(2x - 1)\,\cos t,\nonumber
\end{align}
with $(x, y) \in (0, 1) \times (0, 1)$ along with the Dirichlet boundary condition. Here, the fluid viscosity $\nu$=1, 
time interval $(0,1]$ with final time $T=1$. 
 \end{example}
 \noindent
Tables 1, 2 represent the convergence rates for velocity in $L^{\infty}(\bL^2)$, $L^{\infty}(\bH^1)$-norms, respectively, 
and Table 3 depicts the convergence rates for pressure in $L^{\infty}(L^2 )$-norm for different values of $\kappa$.
The numerical convergence rates presented in tables validate the theoretical findings obtained in Theorem \ref{Tbd1}. 
Moreover, it can be inferred that the numerical results still hold true as $\kappa \rightarrow 0$.\\
\noindent
In Tables 4, 5, we present the velocity error in $L^{\infty}(\bL^2)$-norm and the pressure error in $L^{\infty}(L^2)$-norm, 
respectively, 
for different values of $\kappa$ and fixed $\nu=0.01$. It can be observed from the tables that the 
velocity and pressure errors 
are quite high and are not stable for the mesh size $h=1/2,\,1/4$ for the Navier-Stokes system with $\kappa=0$ and  $\nu=0.01$.
Therefore, more mesh refinement is needed to achieve the desired accuracy. To overcome this issue, we introduce 
a reasonably small presence of $\kappa$ to the Navier-Stokes system and make the system more regularized. Therefore, 
in this case 
by introducing a significantly small value of $\kappa$, $\kappa\in \{10^{-2},\,1\}$ to the Navier-Stokes system, 
the errors of the desired accuracy are achieved at a coarser mesh $h\in \{1/2,\,1/4\}$ with 
much less computational efforts. 
Note that, a significantly small value of $\kappa$ means that here the presence of 
$\kappa=10^{-4}$ in the Navier-Stokes system does not provide the desired accuracy as 
the errors are still quite high and are not stable for the 
mesh size $h=1/2,\,1/4$. The results in tables 4, 5 validate the 
fact that the Kelvin-Voigt model can be thought of as a $\kappa$ regularization 
of the Navier-Stokes model. \\
\begin{example}\label{7e2}
  In this example, we take right-hand side function $\f = 0$, initial condition
$\bu_0 = (10x^2(x -1)^2 y(y -1)(2y -1), -10y^2(y-1)^2 x(x -1)(2x -1), y)$, $\nu = 1$ and
$\kappa = 1$.  Figure 1 represents velocity plots of the Kelvin-Voigt model and the Navier-Stokes model for final time $T = 4$. 
 We observe that the Kelvin-Voigt fluid velocity tends to zero at a slower rate in comparison to the  
 Navier-Stokes fluid velocity. This confirms the fact that after instantaneous removal of the
forces, the velocity of the Kelvin-Voigt fluid does not vanish instantaneously as in the case of the Navier-Stokes fluid.
 \end{example}
\begin{example}
\label{7e3}
This example deals with the benchmark problem lid-driven cavity flow on a unit square with zero body force. 
 Here, no-slip boundary conditions are considered everywhere except non zero velocity
$(u_1, u_2 ) = (1, 0)$ on the upper part of the boundary, that is, the lid of the square is moving with a velocity $(1,\,0)$. 
 For numerical results, we have considered lines $(0.5, y)$ and $(x, 0.5)$, final time $T = 40$, $h = 1/32$ and $\nu = 1$. 
In figure 2, we have compared the Kelvin-Voigt velocity to the steady-state Navier-Stokes velocity for large time 
 and different values of  $\kappa=\{1,0.001,0.00001\}$. The plots depict that the Kelvin-Voigt solution
converges to the steady state solution for large time and as $\kappa\rightarrow 0$.
\end{example}
 
\section{\normalsize \bf Summary}
The article discusses some new higher order regularity estimates for the weak solution 
which are valid for all time $t>0$
and as $\kappa\rightarrow 0$. Semidiscrete optimal error estimates are derived for the velocity in 
$L^{\infty}(\bL^2)$, $L^{\infty}(\bH_0^1)$-norms and for the pressure in $L^{\infty}(\bL^2)$-norm. Further, under 
uniqueness condition, these estimates are shown uniformly in time. Note that, the constants appearing in 
{\it a priori} error bounds are made independent of inverse powers of $\kappa$ by introducing weight functions 
in powers of $t$.  In fact, an introduction 
of these weight functions takes care of regularity issues at time $t=0$. Further, the 
backward Euler method is applied for the complete discretization of the model and fully discrete optimal error estimates are 
derived. Finally, the article is concluded by presenting 
some numerical results which validate our theoretical observations. \\  
\noindent{{\bf Acknowledgement:}  The authors thank Professor Amiya K. Pani for his valuable comments and 
suggestions regarding this.

\small{
\begin{table}[ht!]
 \centering
 \begin{tabular}{|l|l|l|l|l|l|l|l|l|}
 \hline
  $h$   &\small{$\|\bu(t_n)-\bU^n\|_{\bL^2}$}  &\small{$\|\bu(t_n)-\bU^n\|_{\bL^2}$}& \small{$\|\bu(t_n)-\bU^n\|_{\bL^2}$}
        & \small{$\|\bu(t_n)-\bU^n\|_{\bL^2}$}  \\
        &$\kappa=1$                   &$\kappa=10^{-3}$&$\kappa=10^{-6}$&$\kappa=10^{-9}$                            \\
 \hline
 \hline 1/4  & 1.356182 & 1.575487   & 1.575615 & 1.575615 \\
 \hline 1/8  & 1.721053 & 1.784026   & 1.784008 & 1.784008 \\
 \hline 1/16 & 1.879115 & 1.895495   & 1.895489 & 1.895489  \\
 \hline 1/32 & 1.946837 & 1.950956   & 1.950955 & 1.950955  \\
 \hline
 \end{tabular}
 \vspace{.1cm}
 \caption{ Convergence rates for backward Euler method with $k=\mathbb{O}(h^2)$ and $\kappa\rightarrow 0$. }
 \end{table}
 }
 \small{
\begin{table}[ht!]
 \centering
 \begin{tabular}{|l|l|l|l|l|l|l|l|l|}
 \hline
  $h$   &\small{$\|\bu(t_n)-\bU^n\|_{\bH^1}$}  &\small{$\|\bu(t_n)-\bU^n\|_{\bH^1}$}& \small{$\|\bu(t_n)-\bU^n\|_{\bH^1}$}
        & \small{$\|\bu(t_n)-\bU^n\|_{\bH^1}$}  \\
        &$\kappa=1$                   &$\kappa=10^{-3}$&$\kappa=10^{-6}$&$\kappa=10^{-9}$                            \\
 \hline
 \hline 1/4  & 0.666127 & 0.810563   & 0.810686 & 0.810686 \\
 \hline 1/8  & 0.856391 & 0.903887   & 0.903874 & 0.903874 \\
 \hline 1/16 & 0.939755 & 0.952619   & 0.952616 & 0.952616  \\
 \hline 1/32 & 0.973984 & 0.977304   & 0.977303 & 0.977303  \\
 \hline
 \end{tabular}
 \vspace{.1cm}
 \caption{ Convergence rates for backward Euler method with $k=\mathbb{O}(h)$ and $\kappa\rightarrow 0$. }
 \end{table}
 }
 
  \small{
\begin{table}[ht!]
 \centering
 \begin{tabular}{|l|l|l|l|l|l|l|l|l|}
 \hline
  $h$   &\small{$\|p(t_n)-P^n\|_{\bL^2}$}  &\small{$\|p(t_n)-P^n\|_{\bL^2}$}& \small{$\|p(t_n)-P^n\|_{\bL^2}$}
        & \small{$\|p(t_n)-P^n\|_{\bL^2}$}  \\
        &$\kappa=1$                   &$\kappa=10^{-3}$&$\kappa=10^{-6}$&$\kappa=10^{-9}$                            \\
 \hline
 \hline 1/4  & 0.935384 & 0.934915   & 0.934891 & 0.934891 \\
 \hline 1/8  & 0.962732 & 0.960305   & 0.960305 & 0.960305 \\
 \hline 1/16 & 0.982119 & 0.980954   & 0.980955 & 0.980955  \\
 \hline 1/32 & 0.991049 & 0.990312   & 0.990312 & 0.990312  \\
 \hline
 \end{tabular}
 \vspace{.1cm}
 \caption{ Convergence rates for backward Euler method with $k=\mathbb{O}(h^2)$ and $\kappa\rightarrow 0$. }
 \end{table}
 }
  \small{
\begin{table}[ht!]
 \centering
 \begin{tabular}{|l|l|l|l|l|l|l|l|l|}
 \hline
  $h$   &\small{$\|\bu(t_n)-\bU^n\|_{\bL^2}$}  &\small{$\|\bu(t_n)-\bU^n\|_{\bL^2}$}& \small{$\|\bu(t_n)-\bU^n\|_{\bL^2}$}
        & \small{$\|\bu(t_n)-\bU^n\|_{\bL^2}$} \\
        &$\kappa=1$                   &$\kappa=10^{-2}$&$\kappa=10^{-4}$&  $\kappa=0$\\     
 \hline
 \hline 1/2  & 0.156344 & 5.977491   & 29.31149 &   7.060819 \\
 \hline 1/4  & 0.083167 & 3.223425   & 26.16600 & 35.799670\\
 \hline 1/8  & 0.027082 & 1.425267   & 1.810881 &   1.798025  \\
 \hline 1/16 & 0.007577 & 0.432046   & 0.500769 &   0.494285  \\
 \hline 1/32 & 0.001991 & 0.116965   & 0.130588 &   0.128835\\
\hline
 \end{tabular}
 \vspace{.1cm}
 \caption{ Regularization effect on velocity errors for backward Euler method with $\nu=0.01$, $k=\mathbb{O}(h^2)$. }
 \end{table}
 }
 \small{
\begin{table}[ht!]
 \centering
 \begin{tabular}{|l|l|l|l|l|l|l|l|l|l}
 \hline
  $h$   &\small{$\|p(t_n)-P^n\|_{\bL^2}$}  &\small{$\|p(t_n)-P^n\|_{\bL^2}$}& \small{$\|p(t_n)-P^n\|_{\bL^2}$}
        &\small{$\|p(t_n)-P^n\|_{\bL^2}$} \\
       &$\kappa=1$                   &$\kappa=10^{-2}$&$\kappa=10^{-4}$ & $\kappa=0$ \\
 \hline
 \hline 1/2  & 3.057302 & 7.657372   & 160.344805 &   885.464339 \\
 \hline 1/4  & 1.290707 & 3.121213   & 204.368122 & 419.021951 \\
 \hline 1/8  & 0.589169 & 1.386778   &   1.1096204&    1.089742 \\
 \hline 1/16 & 0.280954 & 0.327177   &   0.2832502&   0.283915 \\
 \hline 1/32 & 0.137575 & 0.139234   &   0.1419708 &    0.142181\\
\hline
 \end{tabular}
 \vspace{.1cm}
 \caption{ Regularization effect on pressure errors for backward Euler method with $\nu=0.01$, $k=\mathbb{O}(h^2)$. }
 \end{table}
 }

 \begin{figure}
 		\centering
 	\subfloat[First component of velocity]{\label{sfig:a1}
 		\includegraphics[width=42mm]{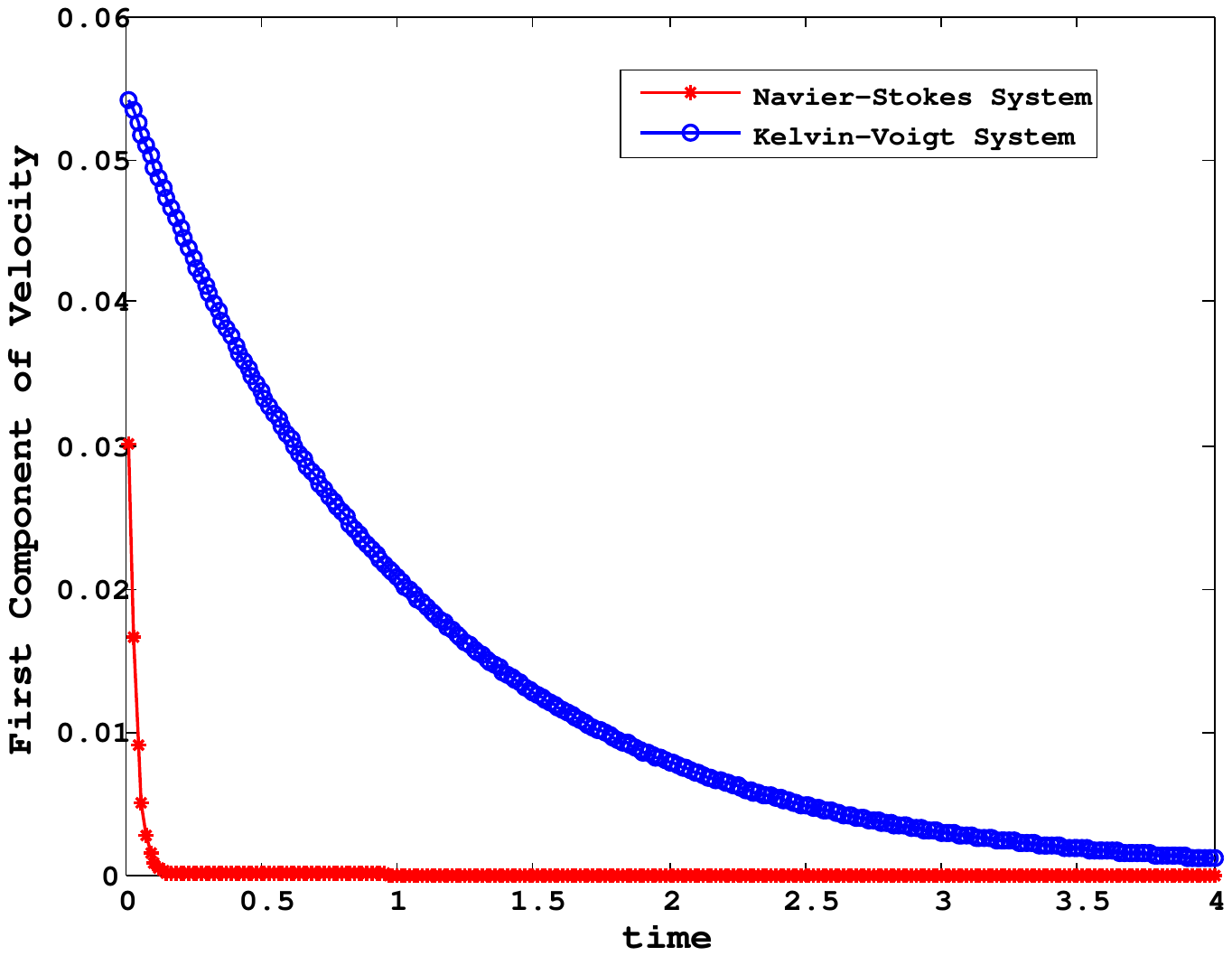}}
 	\hspace{3cm}
 	\subfloat[Second component of velocity]{\label{sfig:b1}\includegraphics[width=42mm]
 		{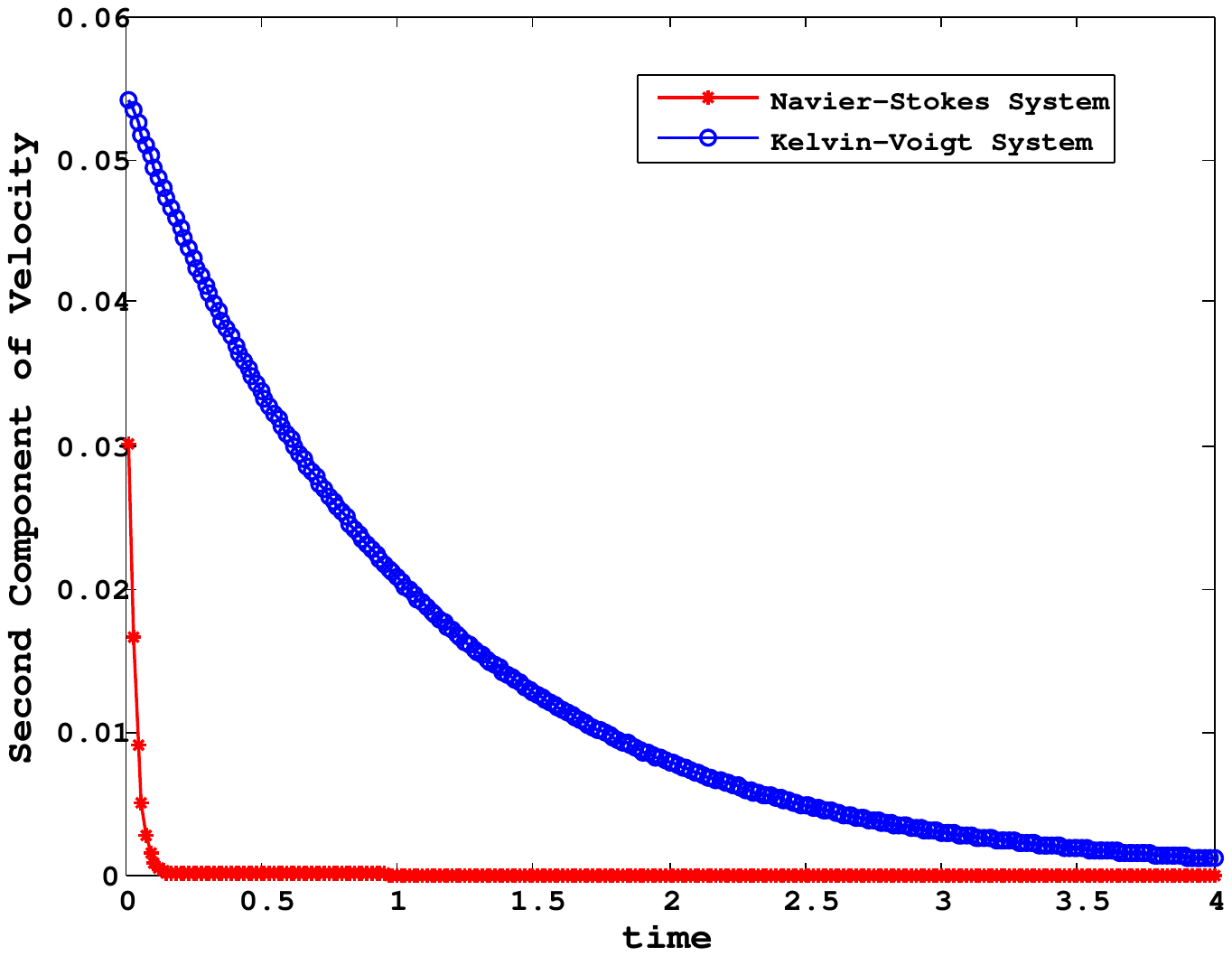}}\
 	\caption{ Comparison of Navier-Stokes velocity and Kelvin-Voigt velocity components for Example \ref{7e2}.}
 	\label{fig1}
	\centering
	\subfloat[First component of velocity for line $x=0.5$]{\label{sfig:a11}
	\includegraphics[width=42mm]{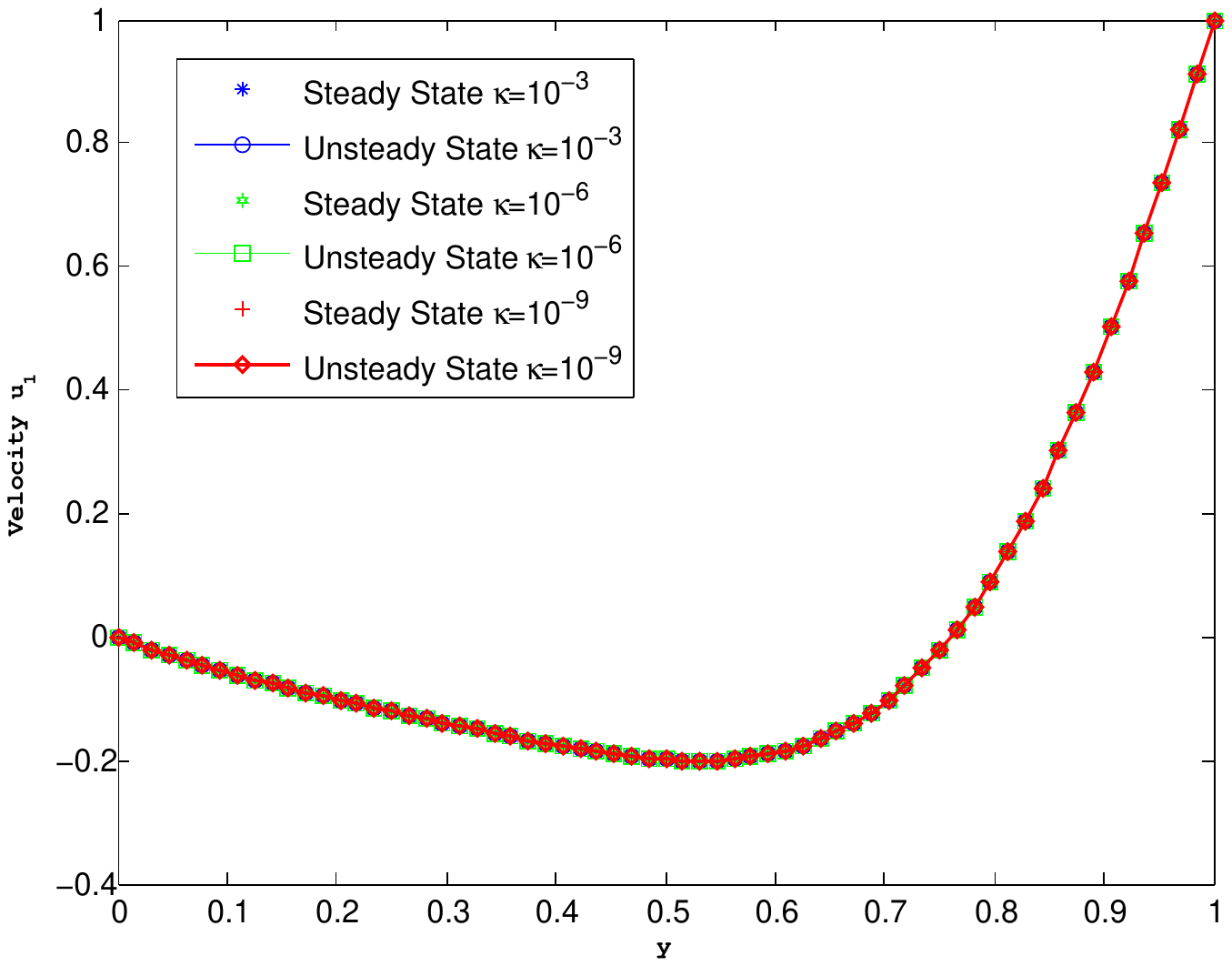}}
	\hspace{4cm}
	\subfloat[Second component of velocity for  line $x=0.5$]{\label{sfig:b11}\includegraphics[width=42mm]
	{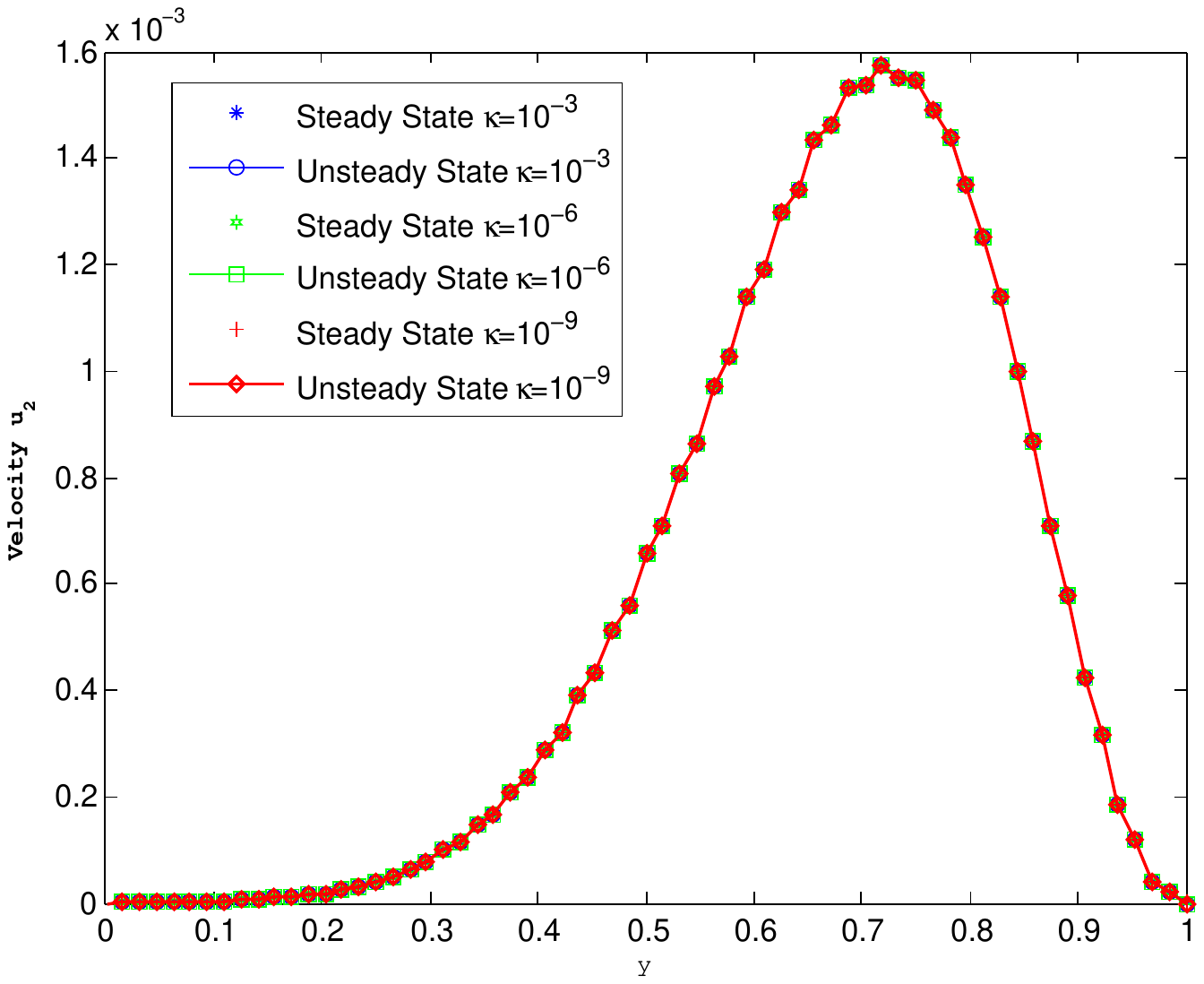}}\\
	\subfloat[First component of velocity for line $y=0.5$]{\label{sfig:c1}\includegraphics[width=42mm]
	{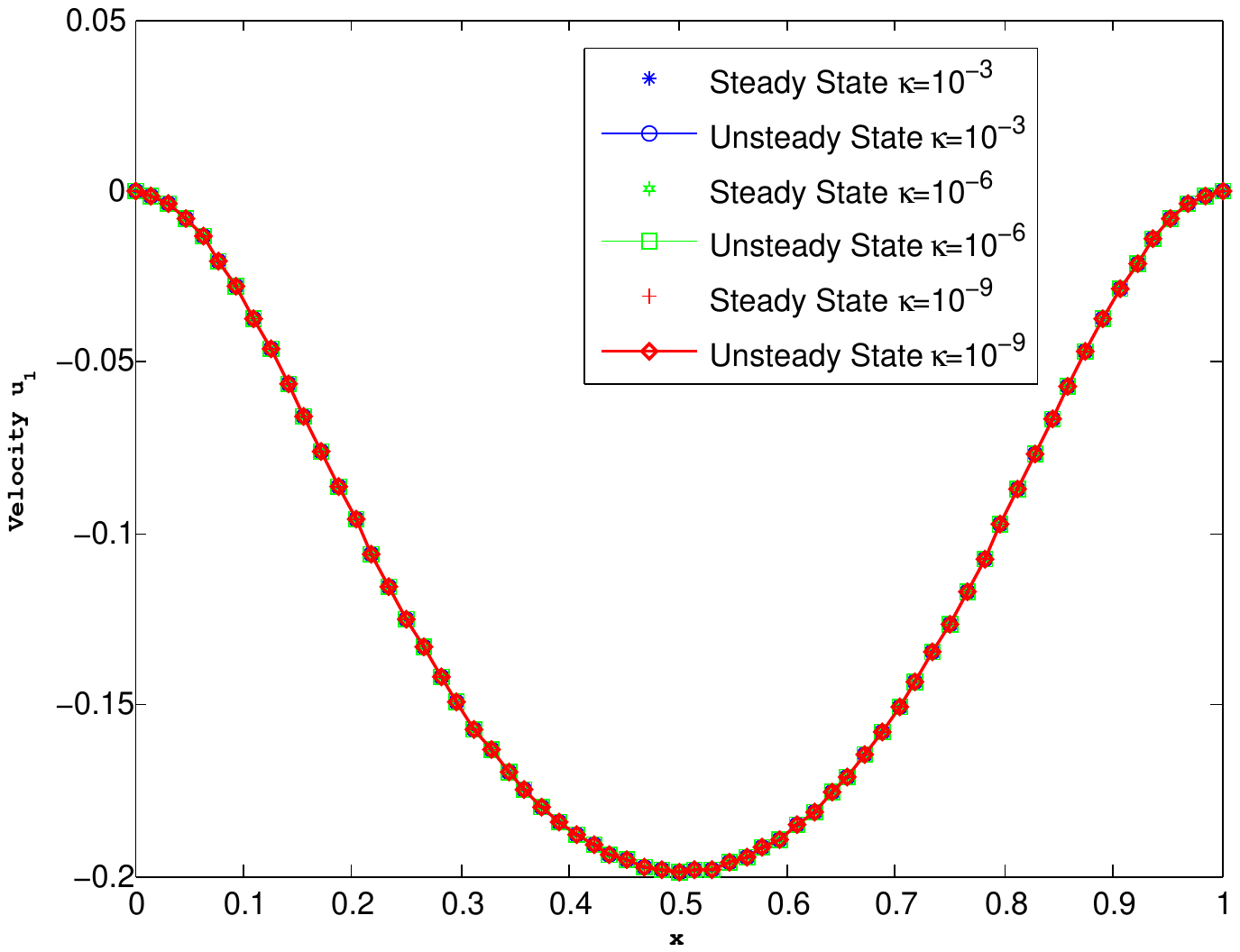}}
	\hspace{4cm}
	\subfloat[Second component of velocity for line $y=0.5$]{\label{sfig:d1}\includegraphics[width=42mm]
	{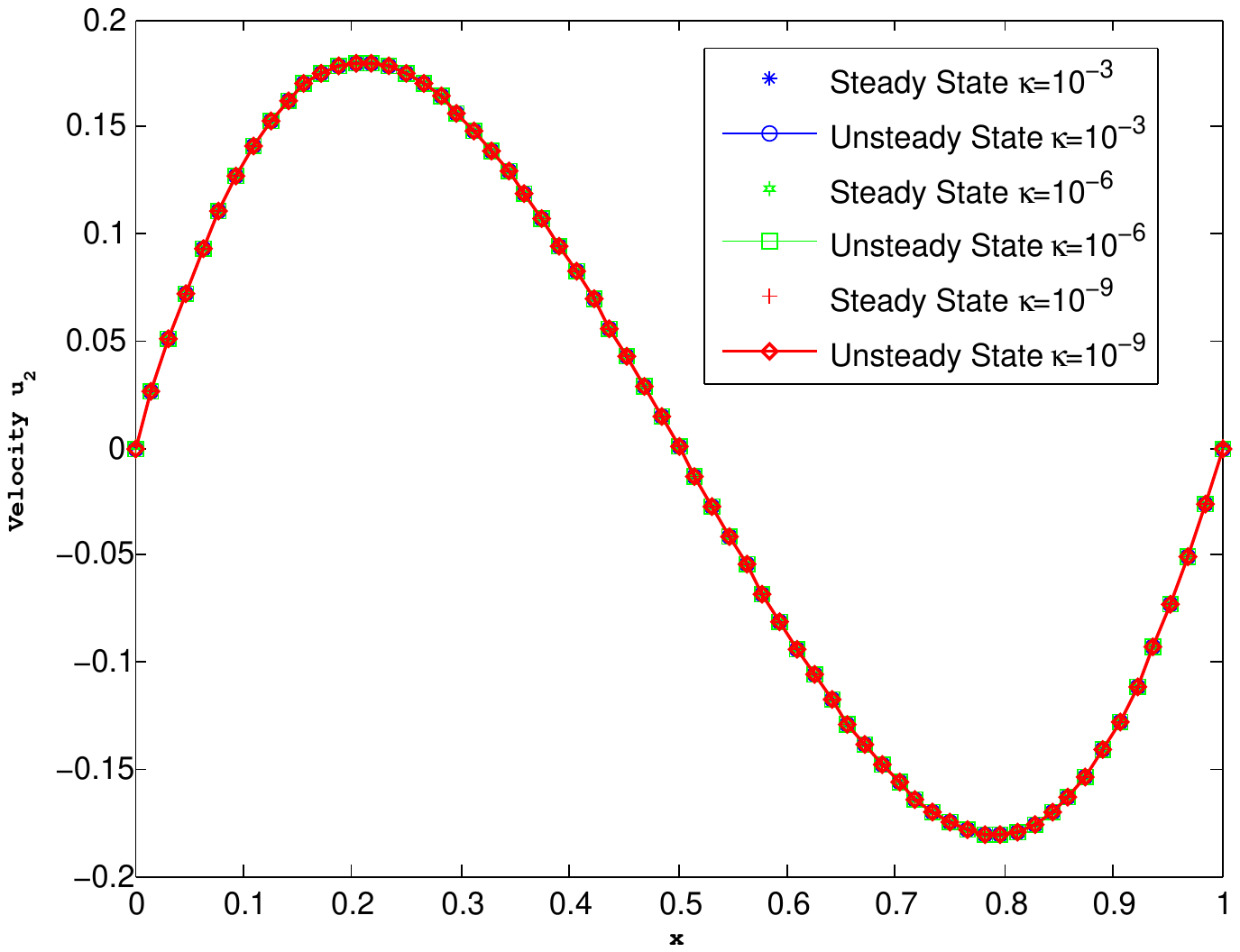}}\
	\caption{ Velocity components for lid-driven cavity flow in Example \ref{7e3}.}
	\label{fig2}
\end{figure}

\end{document}